\documentclass[10pt]{article}

\usepackage{amsmath}
\usepackage{amssymb,amsfonts,bm,bbm,nicefrac}
\usepackage{color,graphicx,caption,subfig}
\usepackage{xspace}
\usepackage{placeins}

\usepackage{enumitem}

\newtheorem{thm}{Theorem}[section]

\newtheorem{remark}[thm]{Remark}

\newcommand{\myinclude}[2]{{\includegraphics{pics/ellipses_a_N10_Ex3}}}

\newcommand\dps{\displaystyle }

\newcommand{\jl}{[\![}
\newcommand{\jr}{]\!]}
\newcommand{\jmp}[1]{\jl#1\jr}

\newcommand{\cd}{\cdot}

\newcommand{\Oo}{\Omega_1}
\newcommand{\Oi}{\Omega_i}
\newcommand{\Oj}{\Omega_j}
\newcommand{\OM}{\Omega_M}

\newcommand{\Oz}{\Omega_0}
\newcommand{\Bz}{\B_R}

\newcommand{\Oinf}{ \Omega_\infty}
\newcommand{\dom}{\Omega}
\newcommand{\B}{\mathsf{B}}

\newcommand{\Gi}{\Gamma_i}
\newcommand{\Gj}{\Gamma_j}
\newcommand{\Ginf}{\Gamma_\infty}
\newcommand{\Gz}{\Gamma_0}

\renewcommand{\ni}{n_i}
\newcommand{\ninf}{n_\infty}
\newcommand{\nz}{n_0}

\newcommand{\IntM}{\mathcal I}

\newcommand{\real}{\mathbb R}
\newcommand{\N}{\mathbb N}
\newcommand{\R}{\mathbb R}
\newcommand{\VNi}{V_{N,i}}
\newcommand{\VNj}{V_{N,j}}
\renewcommand{\L}{\mathsf{L}}
\renewcommand{\H}{\mathsf{H}}
\newcommand{\Kv}{\vec{K}}

\newcommand{\grad}{\nabla}
\renewcommand{\div}{\mbox{div}}
\newcommand{\DtraceOp}[1]{\gamma_{#1}}

\newcommand{\Dtrace}[2]{\gamma_{#1}(#2)}

\newcommand{\Ylm}{\mathcal Y_{\ell m}}
\newcommand{\Ybarlm}{\overline{\mathcal Y}_{\ell m}}
\newcommand{\Yom}{\mathcal Y_{1 m}}
\newcommand{\Yomp}{\mathcal Y_{1 m'}}
\newcommand{\Ylpmp}{\mathcal Y_{\ell ' m'}}
\newcommand{\Ybarlpmp}{\overline{\mathcal Y}_{\ell ' m'}}
\newcommand{\coeff}[1]{[\mathsf{#1}]_\ell^m}

\makeatletter
\newcommand{\pushright}[1]{\ifmeasuring@#1\else\omit\hfill$\displaystyle#1$\fi\ignorespaces}
\newcommand{\pushleft}[1]{\ifmeasuring@#1\else\omit$\displaystyle#1$\hfill\fi\ignorespaces}
\makeatother

\newcommand{\SL}{\widetilde{\mathcal S}}

\newcommand{\SLO}{\mathcal S}
\newcommand{\DLO}{\mathcal D}

\newcommand{\Lop}{\mathcal L}

\newcommand{\cJ}{\mathcal J}
\newcommand{\cJR}{\mathcal J_R}
\newcommand{\cJRN}{{\mathcal J}^{\mathbb{A},N}_R}
\newcommand{\cJRNp}[1]{{\mathcal J}^{\mathbb{A},N}_{R,#1}}

\newcommand{\A}{\mathsf{a}}
\newcommand{\Ainf}{\mathsf{a}_\infty}

\newcommand{\ANR}[1]{\mathsf{a}_{R,N}^{#1}}
\renewcommand{\AA}{\mathsf{A}}
\newcommand{\II}{\mathsf{I}}
\newcommand{\cA}{\mathcal{A}}

\renewcommand{\vec}[1]{{\bm{ #1}}}
\newcommand{\lambdav}{\vec{\lambda}}

\newcommand{\kappav}{\vec{\kappa}}
\newcommand{\tauv}{\vec{\tau}}
\newcommand{\f}{\mathsf{f}}
\newcommand{\fv}{\vec{\mathsf{f}}}
\newcommand{\hv}{\vec{\mathsf{h}}}
\newcommand{\Psiv}{\vec{\Psi}}
\newcommand{\sv}{\vec{\mathsf{s}}}
\newcommand{\p}{\mathsf{p}}
\renewcommand{\u}{\mathsf{u}}
\newcommand{\g}{\mathsf{g}}
\newcommand{\h}{\mathsf{h}}
\renewcommand{\v}{\mathsf{v}}
\newcommand{\I}{\mathsf{I}}
\renewcommand{\S}{\mathsf{S}}

\newcommand{\cha}[1]{\mathbbm{1}_{#1}}
\newcommand{\argmin}[1]{\operatorname*{arg\,min}_{#1}}
\newcommand{\argmax}[1]{\operatorname*{arg\,max}_{#1}}

\title{An embedded corrector problem for homogenization. \\ Part II: Algorithms and discretization}
\author{Eric Canc\`es$^{1,3}$, Virginie Ehrlacher$^{1,3}$, Fr\'ed\'eric Legoll$^{2,3}$, Benjamin Stamm$^4$, Shuyang Xiang$^4$
\\
{\footnotesize $^1$ CERMICS, \'Ecole des Ponts ParisTech, 77455 Marne-La-Vall\'ee Cedex 2, France}
\\
{\footnotesize $^2$ Laboratoire Navier, \'Ecole des Ponts ParisTech, 77455 Marne-La-Vall\'ee Cedex 2, France}
\\
{\footnotesize $^3$ Inria Paris, MATHERIALS project-team, 2 rue Simone Iff, CS 42112, 75589 Paris Cedex 12, France}
\\
{\footnotesize $^4$ MATHCCES, Department of Mathematics, RWTH Aachen University, Schinkelstrasse 2, D-52062 Aachen, Germany}
}
\date{\today}

\begin{document}

\maketitle

\begin{abstract}
This contribution is the numerically oriented companion article of the work~\cite{CELSX}. We focus here on the numerical resolution of the embedded corrector problem introduced in~\cite{notre_cras,CELSX} in the context of homogenization of diffusion equations. Our approach consists in considering a corrector-type problem, posed on the whole space, but with a diffusion matrix which is constant outside some bounded domain. In~\cite{CELSX}, we have shown how to define three approximate homogenized diffusion coefficients on the basis of the embedded corrector problems. We have also proved that these approximations all converge to the exact homogenized coefficients when the size of the bounded domain increases.

We show here that, under the assumption that the diffusion matrix is piecewise constant, the corrector problem to solve can be recast as an integral equation. In case of spherical inclusions with isotropic materials, we explain how to efficiently discretize this integral equation using spherical harmonics, and how to use the fast multipole method (FMM) to compute the resulting matrix-vector products at a cost which scales only linearly with respect to the number of inclusions. Numerical tests illustrate the performance of our approach in various settings.
\end{abstract}


\pagestyle{myheadings}
\thispagestyle{plain} 
\markboth{}{}

\section{Introduction}

Let $\dom \subset \R^d$ be a smooth and bounded domain for some $d\in \N^\star$. We consider the standard elliptic problem 
$$
\left\{
\begin{array}{rcll}
-\hbox{\rm div}\left[ \mathbb A_\varepsilon \, \nabla u_\varepsilon \right] = f & \quad \text{ in $\dom$}, 
\\
\quad u_\varepsilon = 0 & \quad \text{ on $\partial \dom$},
\end{array}
\right.
$$
with a highly oscillatory coefficient of the form $\mathbb{A}_\varepsilon(x)={\mathbb A}(x/\varepsilon)$, where $\varepsilon$ denotes a small-scale parameter and ${\mathbb A}: \R^d \to \R^{d\times d}$ is a field of symmetric positive definite matrices, which is uniformly bounded and coercive. The field ${\mathbb A}$ may be periodic or quasi-periodic. It may also be considered as a realization of a random stationary and ergodic field~\cite{kozlov,papa}. In these latter three cases, it is known from the standard theory of homogenization (see e.g. the textbooks~\cite{Bensoussan,Cioranescu,Jikov} and the review articles~\cite{singapour,engquist-souganidis}) that, as $\varepsilon$ tends to zero, the function $u_\varepsilon$ converges to the solution $u^\star$ to the macroscopic problem
\begin{equation}
\label{eq:div-0}
\left\{
\begin{array}{rcll}
-\mbox{\rm div}\Big( \mathsf A^\star \nabla u^\star \Big) = f &\quad \mbox{ in $\dom$}, 
\\
\quad u^\star = 0 & \quad \text{ on $\partial \dom$},
\end{array}
\right.
\end{equation}
where $\AA^\star \in \R^{d \times d}$ is a constant, positive definite matrix. For practical computations, a so-called corrector problem (which depends on the field ${\mathbb A}$) needs to be solved in order to compute the homogenized matrix~$\AA^\star$, which, in turn, can be used in~\eqref{eq:div-0} for the macroscopic problem without small-scale dependency.

Solving the corrector problem is the key bottleneck in cases beyond the periodic setting. In the quasi-periodic and random stationary cases, one has to solve the corrector problems
\begin{equation}
\label{eq:correctorrandom}
-\mbox{\div}\Big( {\mathbb A}(x)(p + \nabla w(x))\Big) = 0\quad \mbox{ in $\mathcal D'(\R^d)$},
\end{equation}
for $d$ linearly independent vectors $p$ of $\mathbb R^d$ (Problems~\eqref{eq:correctorrandom} are complemented by geometric constraints on $w$, such as $\nabla w$ is quasi-periodic or random stationary, that we do not detail further). Note that the problems~\eqref{eq:correctorrandom} are posed on the whole space $\mathbb R^d$.

In practice, approximate corrector problems, defined on a sequence of increasing truncated domains with appropriate boundary conditions (typically periodic boundary conditions), are often considered to obtain a convergent sequence of approximate homogenized matrices~\cite{Bourgeat,Huet}. Loosely speaking, the larger the size of the truncated domain, the more accurate the corresponding approximation of the homogenized matrix. The use of standard finite element discretizations to tackle these corrector problems may lead to very large discretized problems, the computational cost of which can be prohibitive. 

\medskip

The approach that we have introduced in~\cite{notre_cras,CELSX} is based on a different idea. Rather than truncating the domain on which~\eqref{eq:correctorrandom} is posed, we replace in~\eqref{eq:correctorrandom} the diffusion matrix field ${\mathbb A}$ by the field
$$
\mathcal A^{{\mathbb A},\AA}(x) := \left| 
\begin{array}{ll}
{\mathbb A}(x) & \mbox{ if $x\in \Bz$},\\
\AA &\mbox{ if $x\in \mathbb R^d \setminus \Bz$},
\end{array}\right. 
$$
where $\Bz$ is some bounded domain (typically a sphere of radius $R$) and $\AA$ is some constant matrix to be specified. For some suitable choices of $\AA$, the following corrector-type problem is then considered:
\begin{equation}
\label{eq:pbbase}
-\mbox{\div}\Big(\mathcal A^{{\mathbb A},\AA} (p + \nabla w^{{\mathbb A}, \AA})\Big) = 0 \qquad \mbox{ in $\mathcal D'(\R^d)$}.
\end{equation}
Problem~\eqref{eq:pbbase} is also set on the whole space $\R^d$, but the diffusion matrix is constant outside $\Bz$. In~\cite{notre_cras,CELSX}, we have shown how to define three approximate homogenized matrices from the solution of~\eqref{eq:pbbase}, which all converge to the exact homogenized matrix $\AA^\star$ as $R$ tends to infinity.

\medskip

In this article, we focus on the development of efficient algorithms to solve the embedded corrector problem~\eqref{eq:pbbase} in dimension $d=3$ and to compute the associated effective homogenized coefficients. We consider the case when the following assumptions are satisfied: 
\begin{itemize}
\item at each point, the diffusion matrix ${\mathbb A}$ is isotropic, i.e. proportional to the identity matrix, and therefore reduces to a scalar diffusion coefficient;
\item the diffusion coefficient is piecewise constant and models non-overlapping polydisperse spherical inclusions embedded in a homogeneous material;
\item when $R \to \infty$, the rescaled diffusion coefficients ${\mathbb A}(R \cdot)$ converge in the sense of homogenization to a constant homogenized matrix $\AA^\star$;
\item the homogenized matrix $\AA^\star$ is proportional to identity.
\end{itemize}
In this case, since the diffusion coefficient is piecewise constant, we can reformulate the partial differential equation (PDE)~\eqref{eq:pbbase} in terms of an integral equation whose unknown is a function (the so-called single-layer density) defined on a finite union of spheres (including the sphere $\partial\Bz$). The embedded corrector problems~\eqref{eq:pbbase} are then solved by an algorithm based on a boundary integral formulation~\cite[Chapter~4]{sauter2010boundary}. This integral equation can next be efficiently discretized using spherical harmonics and a (spectral) Galerkin approach. In addition, we resort to an adapted version of the Fast Multipole Method (FMM) to perform fast matrix-vector products at linear cost with respect to the number of spherical inclusions.

Extensions beyond the setting described above will be the subject of further studies. The introduction of the embedded corrector problems~\eqref{eq:pbbase}, and their resolution by an integral equation approach, is motivated by a series of recent works~\cite{Stamm1,Stamm2}, where a very efficient numerical method is proposed to solve Poisson problems with piecewise constant coefficients. We also refer to~\cite{LINDGREN2018712} where a method similar to that presented below has been used, in a completely different context.

\medskip

The approach we discuss here has been initially introduced in~\cite{notre_cras}. The theoretical analysis of the approach, where we show the convergence of our approximations to the exact homogenized matrix, has been performed in~\cite{CELSX}. In this article, we focus on the description of an efficient algorithm to solve~\eqref{eq:pbbase}, and we discuss the numerical performance of our approach.

Note that integral equations have already been introduced in the context of homogenization in~\cite{cazeaux-zahm}. In this article, the authors consider the (standard) corrector problem~\eqref{eq:correctorrandom}, posed on a truncated domain and complemented by periodic boundary conditions. The problem is next recast in a boundary integral formulation (which is different from ours since the PDE of interest is different). This latter formulation is eventually discretized using piecewise constant functions. The theory of hierarchical matrices (namely $\mathcal{H}$-matrices) is used to speed-up the computations.

\medskip

This article is organized as follows. In Section~\ref{sec:homog}, we recall the theoretical results proved in~\cite{CELSX} and our motivation for considering effective homogenized coefficients built upon the embedded corrector problem. We also describe the specific geometrical setting briefly mentioned above, and obtain the PDE formulation~\eqref{eq:PDE_homogen} of the embedded corrector problem in this case. Sections~\ref{sec:IE} and~\ref{sec:disc} are devoted to the presentation of our numerical method. In Section~\ref{sec:IE}, we transform the PDE formulation~\eqref{eq:PDE_homogen} into the equivalent integral equation~\eqref{eq:IntRep} posed on a union of spheres. The discretization of this integral equation by means of spherical harmonics is presented in Section~\ref{sec:disc}. In addition, the fast multipole method (FMM) is used to speed-up the computations, as detailed in Appendix~\ref{ssec:FMM}. Optimization problems have to be solved to obtain our approximate effective coefficients. Section~\ref{sec:coeff} provides algorithmic details on the chosen optimization method. The performance of our numerical method is illustrated on several test cases in Section~\ref{sec:Num}. We see there that our approach provides accurate approximations of the exact homogenized coefficient, and that its cost scales only linearly with respect to the number of inclusions. This allows to numerically consider systems composed of up to $3 \times 10^5$ inclusions using standard laptops, for a computational time of the order of a few tens of minutes.

\section{Motivation: homogenization of isotropic materials with polydisperse spherical inclusions}\label{sec:homog}

\subsection{Theoretical results in the general case}
\label{sec:general_case}

We first recall here some of the theoretical results obtained in the companion article~\cite{CELSX}. For the sake of simplicity, we restrict our presentation to the case of dimension $d=3$.

\medskip

Let $\II$ denote the $3\times 3$ identity matrix, and $0 < \alpha < \beta < +\infty$. Let us denote by 
$$
\mathcal{K}:= \left\{ \AA \in \R^{3\times 3}, \quad \AA \mbox{ symmetric }, \quad \alpha \, \II \leq \AA \leq \beta \, \II \right\},
$$
where the inequalities above have to be understood in the sense of operators. We introduce a measurable field of symmetric positive definite matrices $\mathbb{A}: \R^3 \to \mathcal K$. For any $R>0$, we denote by $\mathbb{A}^R:= \mathbb{A}(R \cdot)$ the rescaled matrix field. Following our companion article (see~\cite[Def.~2.2]{CELSX}), we make the assumption that
\begin{enumerate}[label=({A}\arabic*)]
\item \label{hyp:A1} the family of matrix-valued fields $\left(\mathbb{A}^R\right)_{R>0}$ G-converges to some constant symmetric positive definite matrix $\AA^\star \in \mathcal{K}$ as $R$ goes to infinity. 
\end{enumerate}
We recall that the definition of G-convergence has been introduced by F. Murat and L. Tartar in~\cite{MuratTartar}. In this article, we particularly consider two prototypical situations where Assumption~\ref{hyp:A1} is satisfied: 
\begin{itemize}
 \item [(i)] the field $\mathbb{A}$ is periodic; 
 \item[(ii)] the field $\mathbb{A}$ is a realisation of a random ergodic stationary random field.  
\end{itemize}

The computation of the homogenized matrix $\AA^\star$ is a key challenge in homogenization. In the periodic setting, its value can easily be obtained by solving the so-called corrector problems defined on one periodic unit cell. However, in the stochastic setting, its computation is much more difficult since it requires the resolution of corrector problems defined on the whole space. Only approximations of the homogenized matrix can be computed in practice, and their computation requires in general the resolution of approximate corrector problems defined on truncated domains of large size. 

\medskip

In~\cite{CELSX}, we have proposed three alternative definitions of approximate effective matrices, which we recall hereafter. Let $R>0$ and $\Bz$ denote the centered ball of radius $R$ of $\R^3$. Introduce the spaces
\begin{equation}
\label{eq:def_V0}
V:=\left\{v\in \L^2_{\rm loc}(\R^3), \ \ \nabla v \in \left[\L^2(\R^3)\right]^3\right\}
\qquad \mbox{and} \qquad 
V_0:= \left\{ v \in V, \ \ \int_{\Bz} v = 0\right\}.
\end{equation}
Note that $V_0$ depends on $R$. Since this dependency is irrelevant for what follows, we do not make it explicit in the notation. Let $p\in\R^3$ such that $|p| = 1$ and $\AA_\infty \in \mathcal K$. We consider the unique solution $w^R \in V_0$ to
\begin{equation}\label{eq:DiffEq}
-\div \left(\cA^{\mathbb A,\AA_\infty}(p + \grad w^R) \right) = 0 \quad \mbox{ in $\mathcal{D}'(\R^3)$},
\end{equation}
where 
$$
\cA^{\mathbb A,\AA_\infty}(x):=\left|
\begin{array}{ll}
\mathbb{A}(x) & \mbox{ if $x\in \Bz$}, \\
\AA_\infty & \mbox{ if $x\in \R^3\setminus \Bz$}.
\end{array}
\right.
$$
Problem~\eqref{eq:DiffEq} is called hereafter an {\em embedded corrector problem}. This denomination stems from the fact the diffusion matrix field $\cA^{\mathbb A,\AA_\infty}$ may be seen as the result of the insertion of a ball of radius $R$ of the original material (characterized by the diffusion matrix field $\mathbb{A}$) inside a homogeneous infinite medium characterized by a constant diffusion matrix $\AA_\infty$.

\medskip

The variational formulation of~\eqref{eq:DiffEq} reads as follows: find $w^R \in V_0$ such that
$$
\forall v\in V_0, \quad 
\int_{\Bz} (\nabla v)^\top \mathbb{A} (p + \nabla w^R) 
+ 
\int_{\R^d \setminus \Bz} (\nabla v)^\top \AA_\infty \nabla w^R
- 
\int_{\partial \Bz} \AA_\infty p\cd n \, v = 0,
$$
where $n$ is the outward pointing (w.r.t. $\Bz$) unit normal vector.
Following~\cite{CELSX}, we also introduce the functional $J^{\mathbb{A}}_{R,p}:V_0 \times \mathcal K \to \R$ defined by
\begin{equation}
  \label{eq:Energy}
J^{\mathbb{A}}_{R,p}(v;\AA_\infty)
:=
\frac{1}{|\Bz|} \int_{\Bz} (p+\nabla v)^\top \mathbb A (p+\nabla v) 
+ \frac{1}{|\Bz|} \int_{\Omega_\infty}(\nabla v)^T \AA_\infty \nabla v
- \frac{2}{|\Bz|}\int_{\Gamma_\infty} \AA_\infty p \cd n \, v.
\end{equation}
The function $w^R$ can be equivalently characterized as the unique minimizer of $J^{\mathbb{A}}_{R,p}(\,\cdot\,;\AA_\infty)$:
\[
w^R = \argmin{v\in V_0} J^{\mathbb{A}}_{R,p}(v;\AA_\infty).
\]
Let us also define 
\begin{equation}
  \label{eq:Energy2}
\cJ^{\mathbb{A}}_{R,p}(\AA_\infty) 
:= 
\min_{v\in V_0} J^{\mathbb{A}}_{R,p}(v;\AA_\infty)
=
J^{\mathbb{A}}_{R,p}(w^R;\AA_\infty) \quad \mbox{where $w^R$ is the solution to~\eqref{eq:DiffEq}.}
\end{equation}
Since the application $ p\in \R^3 \mapsto \cJ^{\mathbb{A}}_{R,p}(\AA_\infty)$ is quadratic, there exists a unique symmetric matrix $G^{\mathbb{A}}_R(\AA_\infty) \in \R^{3\times 3}$ such that
$$
\forall p\in \R^3, \qquad p^T G^{\mathbb{A}}_R(\AA_\infty) p = \cJ^{\mathbb{A}}_{R,p}(\AA_\infty).
$$
Let $(e_1,e_2,e_3)$ denote the canonical basis of $\R^3$. In~\cite[Lemma~3.3]{CELSX}, we have proved that the function $\cJR^{\mathbb{A}}: \mathcal K \to \R$ defined by
\begin{equation} \label{eq:somme3}
\forall \AA_\infty \in \mathcal K, \qquad \cJR^{\mathbb{A}}(\AA_\infty) := \frac{1}{3} \sum_{i=1}^3 \cJ^{\mathbb{A}}_{R,e_i}(\AA_\infty) = \frac{1}{3} \mbox{\rm Tr } G^{\mathbb{A}}_R(\AA_\infty)
\end{equation}
is strictly concave and
we have proposed three methods to define approximate homogenized matrices for a finite value of $R$, which we denote by $\AA_R^1$, $\AA_R^2$ and $\AA_R^3$. 

The first method to define an approximate homogenized matrix consists in finding $\AA_R^1\in \mathcal K$ solution to 
\begin{align}
\label{eq:ExCoeff1}
\AA_R^1 = \argmax{\AA_\infty\in \mathcal{K}} \cJR^{\mathbb{A}}(\AA_\infty),
\end{align}
while the second one consists in setting
\begin{align}
\label{eq:ExCoeff2}
\AA_R^2 = G_R^{\mathbb{A}}(\AA_R^1).
\end{align}
Since $\cJR^{\mathbb{A}}$ is strictly concave, $\AA_R^1$ (and thus $\AA_R^2$) is well-defined, and $\AA_R^1$ and $\AA_R^2$ can be jointly obtained. The third method is based on a self-consistent procedure: find $\AA_R^3\in \mathcal{K}$ solution to 
\begin{align}
\label{eq:ExCoeff3}
\AA_R^3 = G_R^{\mathbb{A}}(\AA_R^3).
\end{align}
In general, we have not been able to prove the existence of a matrix $\AA_R^3$ satisfying~\eqref{eq:ExCoeff3}. However, we have been able to establish a weaker existence result (see Section~\ref{sec:isotropic} below) which is sufficient to address the particular case we consider in this article.

In general, the three definitions lead to different values. However, they all converge to the exact homogenized matrix $\AA^\star$ as $R$ goes to infinity (see~\cite[Props.~3.4 and 3.5]{CELSX}): for any $1 \leq i \leq 3$, we have
\[
\lim_{R\to \infty} \AA_R^i = \AA^\star.
\]

\subsection{Isotropic materials with spherical inclusions}
\label{sec:isotropic}

As pointed out above, the embedded corrector problem~\eqref{eq:DiffEq} can, in some cases, be very efficiently solved. This is the situation we consider here. In the sequel, for any $x\in\R^3$ and $r>0$, we denote by $B_r(x)$ the ball of $\R^3$ of radius $r$ centered at $x$. 

\medskip

We make the following additional assumptions on the matrix field $\mathbb{A}$ (see Figure~\ref{fig:Geom}): there exist $\eta>0$, $(x_n)_{n\in \N^\star}\subset \R^3$, $(r_n)_{n\in\N^\star} \subset \R_+^\star$, $(\A_n)_{n\in\N^\star}\subset [\alpha, \beta]$ and $\A_0, \A^\star \in [\alpha, \beta]$ such that
\begin{enumerate} [label=({A}\arabic*)]
  \setcounter{enumi}{1}
\item \label{hyp:A2} for all $n\neq m\in \N^\star$, $\mbox{\rm dist}(B_{r_n}(x_n), B_{r_m}(x_m))\geq \eta$; 
\item \label{hyp:A3} for all $n\in \N^\star$, $\mathbb{A}(x) = \A_n \, \II$ when $x\in B_{r_n}(x_n)$; 
\item \label{hyp:A4} $\mathbb{A}(x) = \A_0 \, \II$ on $\Bz \setminus \bigcup_{n\in\N^\star} B_{r_n}(x_n)$.
\item \label{hyp:A5} $\AA^\star = \A^\star \, \II$. 
\end{enumerate}
In other words, we focus here on the case when the matrix-valued field $\mathbb{A}$ models a material composed only of isotropic phases ($\mathbb{A}$ is everywhere proportional to the identity matrix $\II$), with spherical inclusions embedded into a homogeneous material, and such that the associated homogenized material is also isotropic. The geometry of some realistic heterogeneous materials, such as cement foams or cement-based materials for instance, may be modeled (at least approximately) by spherical inclusions embedded in a homogeneous material (see e.g.~\cite[Figure 9]{bisschop} and~\cite{manila}).

\medskip

Assume now that we are interested in computing the homogenized coefficient $\A^\star$ associated to a matrix-valued field $\mathbb{A}$ satisfying Assumptions~\ref{hyp:A1}--\ref{hyp:A5}. Then, following the results of~\cite{CELSX}, for each value of $R>0$, one can define three approximate effective coefficients $\A_R^1$, $\A_R^2$ and $\A_R^3$, which are scalar versions of~\eqref{eq:ExCoeff1}, \eqref{eq:ExCoeff2} and~\eqref{eq:ExCoeff3}, as follows:
\begin{align}
\label{eq:ExCoeffScal1}
\A_R^1 &= \argmax{\A_\infty\in [\alpha, \beta]} \cJR^{\mathbb{A}}(\A_\infty \, \II),
\\
\label{eq:ExCoeffScal2}
\A_R^2 &= \cJR^{\mathbb{A}}(\A_R^1 \, \II),
\\
\label{eq:ExCoeffScal3}
\A_R^3 \in [\alpha,\beta] \ \ & \text{such that} \ \ \A_R^3 = \cJR^{\mathbb{A}}(\A_R^3 \, \II).
\end{align}
Since $\cJR^{\mathbb{A}}$ is strictly concave, and in view of~\cite[Prop.~3.7]{CELSX}, the above three approximations are well-defined. An easy adaptation of the arguments presented in~\cite{CELSX} yields that, for any $1 \leq i \leq 3$,
$$
\mathop{\lim}_{R\to +\infty} \A_R^i = \A^\star.
$$

The motivation for considering such effective approximations is the following: computing the solution of the embedded corrector problems~\eqref{eq:DiffEq} with $\AA_\infty = \A_\infty \, \II$ for some $\A_\infty>0$ can be done very efficiently, provided that the sphere $\partial \Bz$ does not intersect any of the spherical inclusions $B_{r_n}(x_n)$ for $n\in \N^\star$. Under this assumption, we recast~\eqref{eq:DiffEq} as an interface problem (see~\eqref{eq:PDE_homogen} below), and propose in this article a very efficient numerical method for the resolution of~\eqref{eq:PDE_homogen} and the computation of $\A_R^1$, $\A_R^2$ and $\A_R^3$ (see Sections~\ref{sec:IE}, \ref{sec:disc} and~\ref{sec:coeff} below for details).

\bigskip

\begin{figure}[t!]
\centering
\includegraphics[scale=0.6]{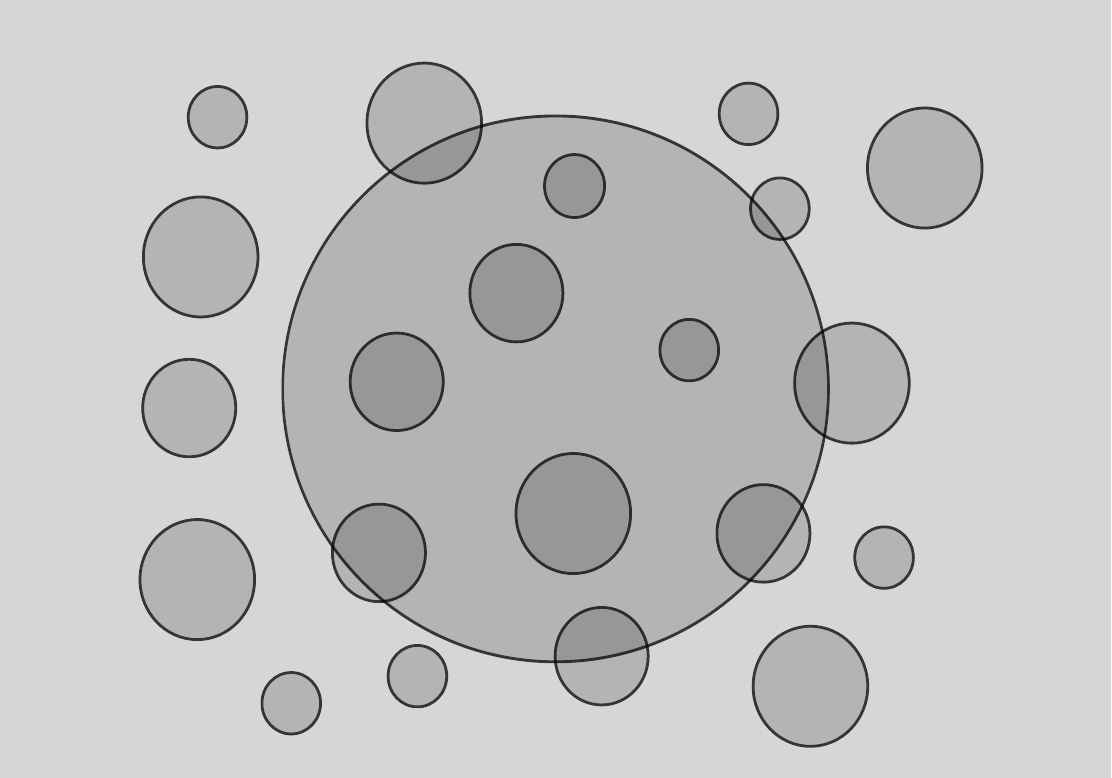}
\includegraphics[scale=0.6]{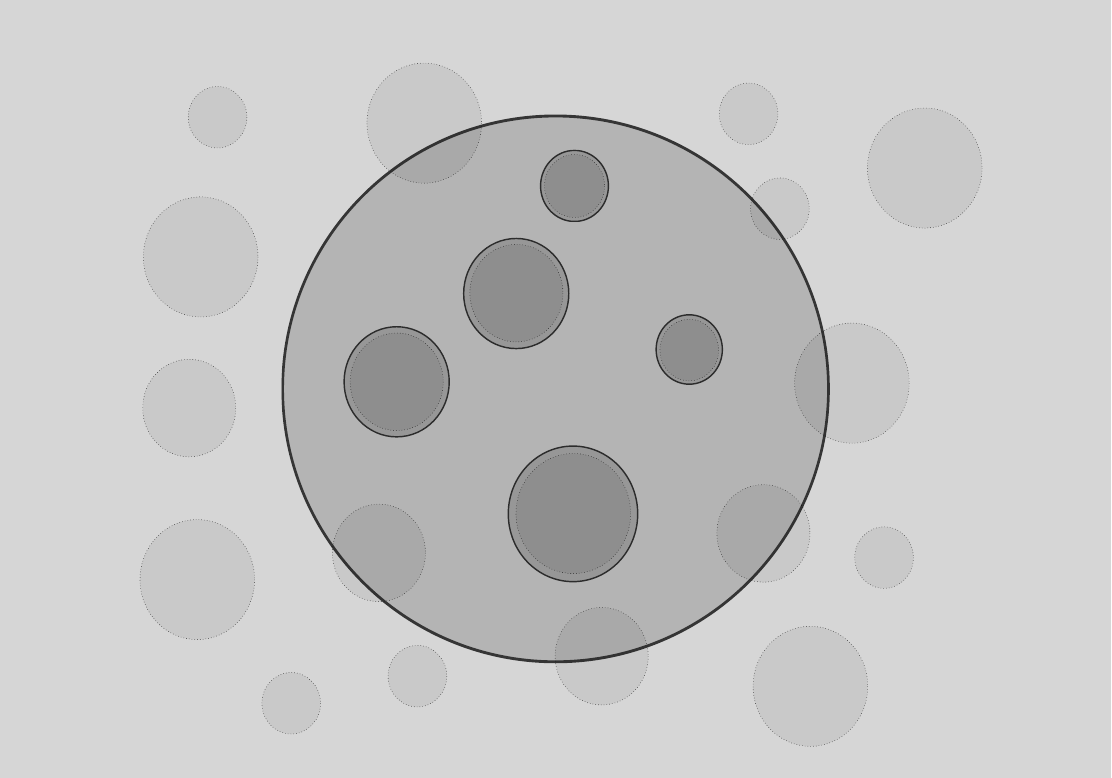}
\caption{Process of restricting and scaling all spherical inclusions inside $\Bz$.}
\label{fig:Geom12}
\end{figure}

Of course, if the matrix-valued field $\mathbb{A}$ satisfies Assumptions~\ref{hyp:A1}--\ref{hyp:A5}, it is not always possible to find arbitrarily large values of $R>0$ such that $\partial \Bz$ does not intersect any of the spherical inclusions of the material (see the left side of Figure~\ref{fig:Geom12} for an illustration). Adapting the algorithm presented in this article to the case when spherical inclusions can intersect with each other and/or when spherical inclusions can intersect with the outer sphere $\partial \Bz$ will be the subject of a future work.

Here, we propose a heuristic procedure which consists, for a given value of $R>0$, in replacing the value of the material coefficient field $\mathbb{A}$ inside the ball $\Bz$ by a modified coefficient field $\overline{\mathbb{A}}$ defined as follows. Let us assume that $\mbox{\rm Card}\{ n\in \N^\star, \ \ B_{r_n}(x_n) \subset \Bz\} = M$ for some $M\in \N^\star$: there are exactly $M$ spherical inclusions that are contained in the ball $\Bz$. Up to reordering the elements of the sequence $(x_n,r_n,\A_n)_{n\in\N^\star}$, we can assume that $\{ n\in \N^\star, \ \ B_{r_n}(x_n) \subset \Bz\} = \{1, \ldots, M\}$ without loss of generality. Let us assume in addition that there are $\widehat{M} \in \N^\star$ balls that intersect with $\Bz$ but do not lie entirely in $\Bz$: 
$$
\mbox{\rm Card} \Big\{ n\in \N^\star, \ \ B_{r_n}(x_n) \cap \Bz \neq \emptyset \ \ \mbox{and} \ \ B_{r_n}(x_n) \not\subset \Bz \Big\} = \widehat{M}.
$$
We denote by $\widehat{x}_1, \ldots, \widehat{x}_{\widehat{M}}$ (respectively $\widehat{r}_1, \ldots, \widehat{r}_{\widehat{M}}$ and $\widehat{a}_1, \ldots, \widehat{a}_{\widehat{M}}$) their centers (respectively their radii and diffusion coefficients).

We define 
\begin{equation}
\label{eq:gamma}
\gamma:= \frac{\sum_{i=1}^M |B_{r_i}(x_i)| + \sum_{j=1}^{\widehat{M}} |\Bz \cap B_{\widehat{r}_j}(\widehat{x}_j)|} {\sum_{i=1}^M |B_{r_i}(x_i)|}.
\end{equation}
The material coefficient $\mathbb{A}$ inside the ball $\Bz$ is then replaced by the modified material coefficient
\begin{equation}\label{eq:scaling}
\overline{\mathbb{A}}(x):= \left|
\begin{array}{ll}
 \A_i \, \II & \mbox{ if $x \in B_{\gamma r_i}(x_i)$ for some $1 \leq i \leq M$},\\
 \A_0 \, \II & \mbox{ if $x\in \Bz \setminus \bigcup_{i=1}^M B_{\gamma r_i} (x_i)$}.
\end{array} \right.
\end{equation}
In other words, in the proposed procedure, the spherical inclusions which intersect with $\partial \Bz$ are deleted and the ones included in $\Bz$ are rescaled by the factor $\gamma \geq 1$. A sketch of this procedure is shown on the right side of Figure~\ref{fig:Geom12}.

Since the spherical inclusions included in $\Bz$ grow in the rescaling process, some of them may no longer be included in $\Bz$ after rescaling, or may intersect with another inclusion. However, since the scaling factor $\gamma$ converges to one as $R$ goes to infinity with rate $\mathcal{O}\left(R^{-2/3}\right)$, we do not observe this problem in practice in our numerical tests, for large values of $R$. Effective coefficients $\A_R^1$, $\A_R^2$ and $\A_R^3$ are then computed using formulas~\eqref{eq:ExCoeffScal1}, \eqref{eq:ExCoeffScal2} and~\eqref{eq:ExCoeffScal3}, using $\overline{\mathbb{A}}$ instead of $\mathbb{A}$. 

\subsection{Geometrical setup and governing equations}\label{sec:prob}

In this section, we describe the geometrical setup we consider and detail the governing equations we solve. The notation introduced here is used in all the sequel.

\medskip

\begin{figure}[t!]
\centering
\def\svgwidth{0.8\textwidth}
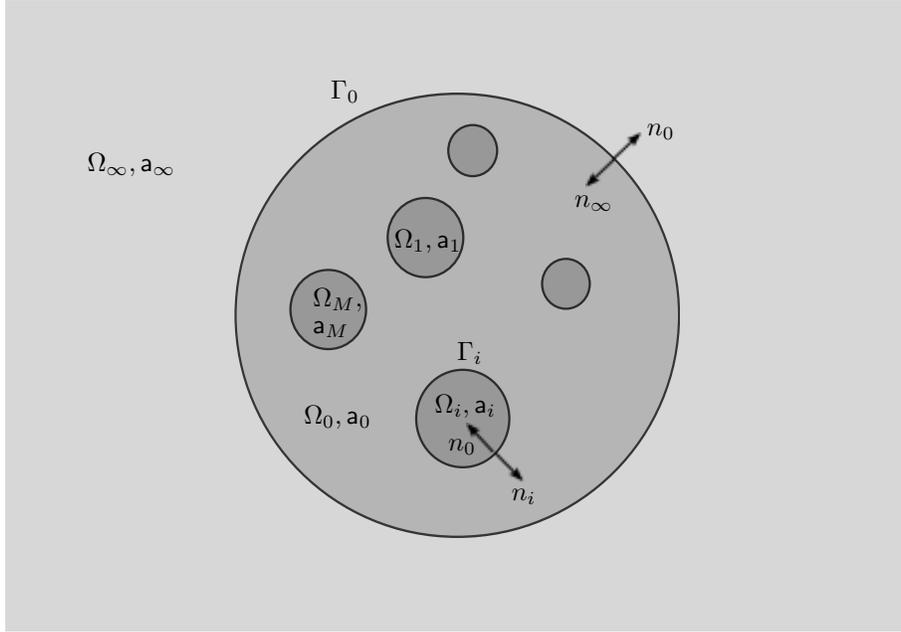
\caption{Geometrical configuration.}
\label{fig:Geom}
\end{figure}

We fix $R>0$, $p\in \R^3$ such that $|p|=1$ and $\A_\infty>0$. Let $\Oo,\ldots,\OM$ be a collection of $M$ non-overlapping balls in $\mathbb R^3$, where each $\Omega_i$ has radius $r_i$ and is centered at $x_i\in\mathbb R^3$ (see Figure~\ref{fig:Geom}). Let $\Bz$ be the ball centered at the origin of radius $R$. As explained above, we assume that $\Omega_i \subset \Bz$ for any $1 \leq i \leq M$. We assume that there exists $\eta >0$ such that 
$$
\mathop{\min}_{1 \leq n\neq m \leq M} \mbox{\rm dist}\left(\Omega_n, \Omega_m \right)\geq \eta \quad \mbox{ and } \quad \mathop{\min}_{1 \leq m \leq M} \mbox{\rm dist}\left(\partial \Bz, \Omega_m \right)\geq \eta.
$$
Set $r_\infty = R$ and $x_\infty=0$ (the origin). Let $\Oinf = \real^3\setminus\Bz$ and $\Oz=\Bz\setminus \cup_{i=1}^M \Oi$. Set $\Gz=\partial\Oz$ and $\Gamma_i=\partial\Omega_i$ for any $i\in\IntM$, with
$$
\IntM := \{1,\ldots,M,\infty\}.
$$
Note that
\begin{align*}
\R^3\setminus \overline{\Oz} &= \bigcup_{i\in\IntM} \Omega_i = \Omega_1\cup \ldots \cup \Omega_M\cup \Oinf,
\\
\Gamma_0 &= \bigcup_{i\in\IntM} \Gamma_i = \Gamma_1\cup \ldots \cup \Gamma_M\cup \Gamma_\infty.
\end{align*}

For all $i\in\IntM$, $n_i$ is the outward pointing (w.r.t. $\Omega_i$) unit normal vector on $\Gamma_i$ and $n_0 = - n_i$. In addition, for all $i \in \IntM \cup \{0\}$, we denote by $\DtraceOp{i}:\mathrm{H}^1(\Omega_i) \to \mathrm{H}^{\nicefrac12}(\Gamma_i)$ the Dirichlet trace operator. For any function $v\in H^1(\Omega_i)$ such that $\Delta v = 0$ in $\Omega_i$, we denote by $\nabla v \cdot n_i$ its (interior with respect to $\Omega_i$) co-normal derivative on $\Gamma_i$ (with outward pointing normal vector).

\medskip

We denote the uniform diffusion parameter inside $\Omega_i$ by $\A_i>0$, for $i\in\IntM \cup \{0\}$. We also define the piecewise-constant diffusion coefficients $\mathbb A: \Bz \to \R$ and $\mathcal A^{{\mathbb A},\A_\infty}: \R^3 \to \R$ as follows:
\begin{align*}
\mathbb A(x) & := \A_0 + \sum_{i=1}^M (\A_i-\A_0) \cha{\Omega_i}(x), \qquad x \in \Bz, \\
\mathcal A^{\mathbb A,\A_\infty}(x) & := \A_0 + (\A_\infty-\A_0) \cha{\Omega_\infty}(x) + \sum_{i=1}^M (\A_i-\A_0) \cha{\Omega_i}(x), \qquad x\in \real^3,
\end{align*}
where $\cha{\Omega_i} $ is the characteristic function of $\Omega_i$ (see Figure~\ref{fig:Geom}). Note that, by definition, 
$$
\mathcal A^{\mathbb A,\A_\infty}(x) = \left|
\begin{array}{ll}
\mathbb A(x) & \mbox{ if $x \in \Bz$},\\
\A_\infty & \mbox{ if $x \in \R^3 \setminus \Bz$}.   
\end{array} \right.
$$

\bigskip

As pointed out above, we focus on the following embedded corrector problem (from now on, we do not make explicit the dependency of $w$ upon $R$): find $w\in V_0$ such that
\begin{equation}\label{eq:DiffEq2}
-\div \left( \cA^{\mathbb A,\A_\infty}(p + \grad w) \right) = 0 \quad \mbox{ in $\mathcal{D}'(\R^3)$}, 
\end{equation}
where the space $V_0$ is defined by~\eqref{eq:def_V0}. Following~\cite{CELSX}, the variational formulation of~\eqref{eq:DiffEq2} reads as: find $w\in V_0$ such that
\begin{equation}
\label{eq:VarForm2}
\forall v\in V_0, \quad 
\int_{\Bz} (\nabla v)^\top \mathbb{A} (p + \nabla w) 
+
\Ainf \int_{\Omega_\infty} (\nabla v)^\top \nabla w
- \Ainf \int_{\Gamma_\infty} p\cd n \, v = 0.
\end{equation}
Following Section~\ref{sec:general_case}, we are interested in computing the quantity (see~\eqref{eq:Energy}--\eqref{eq:Energy2})
$$
\cJ^{\mathbb{A}}_{R,p}(\Ainf) :=\frac{1}{|\Bz|} \sum_{i=0}^M
\A_i\int_{\Oi} |p+\nabla w|^2
+  \frac{\Ainf }{|\Bz|} \int_{\Oinf} |\nabla w|^2 
+ 2 \frac{\Ainf }{|\Bz|}\int_{\Ginf}  p\cd \ninf \, w.
$$
Note that, in the last term, we use the normal vector $\ninf = -n$. Using~\eqref{eq:VarForm2}, it holds that
\begin{align}
\cJ^{\mathbb{A}}_{R,p}(\Ainf) &= \frac{1}{|\Bz|} \sum_{i=0}^M
\left[
\A_i \int_{\Oi}|p|^2 + \A_i\int_{\Oi} p^\top\nabla w
\right]
+ \frac{\Ainf}{|\Bz|} \int_{\Ginf} p\cd \ninf \, w
\nonumber
\\
&= \frac{1}{|\Bz|} \left[\sum_{i=0}^M \A_i |\Oi| + \int_{\Gz} \jmp{\cA^{\mathbb A,\A_\infty}  p \cd n} w\right],
\label{eq:lundi}
\end{align}
where, for all $i\in \IntM$, the jump $\jmp{\cA^{\mathbb A,\A_\infty} p \cd n}$ is defined on $\Gamma_i$ by
$$
\jmp{\cA^{\mathbb A,\A_\infty} p \cd n} := \Dtrace{0}{\cA^{\mathbb A,\A_\infty} p} \cd n_0 + \Dtrace{i}{\cA^{\mathbb A,\A_\infty} p} \cd n_i= (\A_0 \, n_0 + \A_i \, n_i)\cd p.
$$
Thus, the computation of the quantity $ \cJ^{\mathbb{A}}_{R,p}(\Ainf)$ only requires the knowledge of the function $w$ on $\Gamma_0$, and not necessarily in the whole space $\R^3$. This remark motivates the introduction of an integral formulation of Problem~\eqref{eq:DiffEq2}, whose unknown function is the trace of the function $w$ on $\Gamma_0$. Such a formulation is likely to be more efficient from a computational point of view, since it only requires the knowledge of the values of $w$ on the spheres $\Gamma_i$ for $i\in \IntM$, and not in the whole space $\R^3$.

\medskip

In order to use an integral equation formulation, we recast the problem~\eqref{eq:DiffEq2} as an equivalent interface problem: find $w\in V_0$ such that
\begin{equation}\label{eq:PDE_homogen}
  \left\{
  \begin{array}{rll}
\Delta w &= 0
&
\mbox{in $\real^3\setminus \Gamma_0$},
\\
\jmp{w} &= 0 
&
\mbox{on $\Gamma_0$},
\\
\jmp{\cA^{\mathbb A,\A_\infty}\grad w \cd n} &= -\jmp{\cA^{\mathbb A,\A_\infty} p \cd n}  
&
\mbox{on $\Gamma_0$},
\end{array}
\right.
\end{equation}
where, for all $i \in \IntM$, the jumps $\jmp{w}$ and $\jmp{\cA^{\mathbb A,\A_\infty} \grad w \cd n}$
are defined on $\Gamma_i$ by
\begin{align*}
\jmp{w}
&:= 
\Dtrace{i}{w} - \Dtrace{0}{w},
\\
\jmp{\cA^{\mathbb A,\A_\infty} \grad w \cd n}
&:= 
\A_0 \nabla w \cdot n_0  + \A_i \nabla w \cdot n_i.
\end{align*}
In the sequel, we use the notation $\lambda := \gamma_0(w)\in \mathsf{H}^{\nicefrac12}(\Gamma_0)$ for the trace of $w$ on $\Gamma_0$. We also denote by $\lambda_i := \lambda|_{\Gamma_i} = \gamma_i(w) \in \mathsf{H}^{\nicefrac12}(\Gi)$ the trace of $w$ on the sphere $\Gi \subset \Gamma_0$, for any $i\in \IntM$. The equality~\eqref{eq:lundi} thus writes
\begin{equation} \label{eq:nrj_lambda}
\cJ^{\mathbb{A}}_{R,p}(\Ainf) = \frac{1}{|\Bz|} \left[\sum_{i=0}^M \A_i |\Oi| + \int_{\Gz} \jmp{\cA^{\mathbb A,\A_\infty} p \cd n} \lambda\right].
\end{equation}

\section{An integral representation}
\label{sec:IE}

In this section, we derive an integral equation formulation of Problem~\eqref{eq:PDE_homogen}. The unknown function of this integral equation formulation is $\lambda$, the trace of $w$ on $\Gamma_0$. We first start by recalling some basic tools of potential theory in Sections~\ref{sec:loc} and~\ref{sec:glob} that are needed in the derivation of the integral equation. We next obtain in Section~\ref{sec:int_eq} the integral equation~\eqref{eq:IntRep}, which is equivalent to~\eqref{eq:PDE_homogen} and much easier to numerically solve.

\subsection{Local representation}\label{sec:loc}

At a local level, i.e. for each sphere $\Gi$ with $i\in \IntM$, we introduce the local single layer potential operator $\SL_i:\mathsf{H}^{-\nicefrac12}(\Gi)\to \mathsf{H}^1(\R^3\setminus \Gi)$ defined by (see e.g.~\cite{wendland})
\begin{align*}
\forall \sigma_i \in \mathsf{H}^{-\nicefrac12}(\Gi), \quad \forall x \in \R^3 \setminus \Gi, \qquad (\SL_i\sigma_i)(x) := \int_{\Gi} \frac{\sigma_i(s)}{|x-s|} \, ds.
\end{align*}
The single layer potential $\SL_i\sigma_i$ is harmonic in $\R^3\setminus \Gi$, vanishes at infinity and is continuous across the interface $\Gi$. We also recall that, up to a multiplicative constant, $\sigma_i$ is the jump of the normal derivative of $\SL_i\sigma_i$ across $\Gi$:
$$
\sigma_i = \frac{1}{4 \pi} \left( \nabla \left(\SL_i\sigma_i\right) \cdot n_i + \nabla \left(\SL_i\sigma_i\right) \cdot n_0 \right).
$$
Since $\SL_i\sigma_i$ is continuous across $\Gi$, we can define the local single layer boundary operator $\SLO_i:\mathsf{H}^{-\nicefrac12}(\Gi)\to \mathsf{H}^{\nicefrac12}(\Gi)$ by restricting the single layer potential to $\Gi$: we set $\SLO_i\sigma_i := \left. \SL_i\sigma_i \right|_{\Gi}$, that is
\begin{align*}
\forall s \in \Gi,
\qquad
(\SLO_i\sigma_i)(s) 
=  
\Dtrace{i}{\SL_i\sigma_i}(s)
=  
\Dtrace{0}{\SL_i\sigma_i}(s)
=
\int_{\Gi} \frac{\sigma_i(s')}{|s-s'|} \, ds'.
\end{align*}
We recall that the local single layer boundary operator $\SLO_i$ is invertible~\cite[Chapter~4]{sauter2010boundary}.

\medskip

For any $1\leq i \leq M$, denote by $v_i$ the harmonic extension of $\lambda_i$ inside $\Omega_i$. Of course, in view of~\eqref{eq:PDE_homogen}, it holds that $v_i = w|_{\Omega_i}$. Moreover, the function $v_i$ can be represented by a single layer potential, i.e. there exists a unique $\sigma_i \in \mathsf{H}^{-\nicefrac12}(\Gi)$ such that
$$
w = v_i=\SL_i\sigma_i \qquad \mbox{ in $\Omega_i$},
$$
and $\sigma_i$ satisfies the integral equation
\begin{equation}\label{eq:Sigi}
\SLO_i \sigma_i = \lambda_i \qquad \mbox{on $\Gamma_i$}.
\end{equation}
Note that $\SL_i \sigma_i$ coincides with $w$ in $\overline{\Oi}$ but not in $\R^3 \setminus \overline{\Oi}$ ($w$ is not harmonic in $\R^3 \setminus \overline{\Oi}$). Furthermore, the interior (with respect to $\Oi$) co-normal derivative of $v_i$ (with outward pointing normal vector) on $\Gi$ is given by
\begin{equation}\label{eq:NeuTrInt0}
\nabla w\cd\ni = \nabla v_i\cd \ni = (2\pi + \DLO_i^\star)\sigma_i \quad \text{on $\Gi$},
\end{equation}
where $\DLO_i^\star$ is the adjoint of the double layer boundary operator which, in the special case of spherical domains, is self-adjoint and satisfies
$$
\DLO_i^\star = \DLO_i = -\frac{1}{2r_i} \SLO_i
$$
for a sphere of radius $r_i$. This last equation, together with~\eqref{eq:Sigi} and~\eqref{eq:NeuTrInt0}, yields
\begin{equation}\label{eq:NeuTrInt}
\nabla w\cd\ni = \nabla v_i \cd n_i = \left(2\pi \SLO_i^{-1} - \frac{1}{2r_i} \right) \lambda_i.
\end{equation}
The situation is slightly different for the region $\Oinf$. Indeed, let $v_\infty$ denote the harmonic extension of $\lambda_\infty$ inside $\Omega_\infty = \R^3 \setminus \Bz$ which decays to $0$ at infinity. Of course, $v_\infty = w|_{\Oinf}$, and there exists a unique $\sigma_\infty \in \mathsf{H}^{-\nicefrac12}(\Gamma_\infty)$ such that
$$
w = v_\infty = \SL_\infty  \sigma_\infty \qquad \mbox{ in $\Omega_\infty$}.
$$
The interior (with respect to $\Oinf$ and thus exterior with respect to $\Oz$) co-normal derivative of $v_\infty$ on $\Ginf$ is then given by
\begin{equation}
\label{eq:NeuTrExt}
\nabla w\cd n_\infty =  \grad v_\infty \cd n_\infty
= (2\pi - \DLO_\infty^\star)\sigma_\infty
= \left(2\pi + \frac{1}{2R} \SLO_\infty\right)\sigma_\infty
= \left(2\pi \SLO_\infty^{-1} + \frac{1}{2R} \right)\lambda_\infty.
\end{equation}

\subsection{Global representation}\label{sec:glob}

We now consider the region $\Omega_0$. The global single layer potential $\SL_{\tt G}:\mathsf{H}^{-\nicefrac12}(\Gz) \to \mathsf{H}^1(\R^3\setminus \Gz)$ is defined as follows: for any $\nu \in \mathsf{H}^{-\nicefrac12}(\Gz)$, 
\begin{equation} \label{eq:def_tildeSG}
\forall x \in \R^3 \setminus \Gz, \qquad (\SL_{\tt G} \nu)(x) := \sum_{i\in\IntM} (\SL_i \nu_i)(x),
\end{equation}
where, for all $i\in \IntM$, $\nu_i = \nu|_{\Gi}$ on each $\Gi$. We can also consider the global single layer boundary operator $\SLO_{\tt G} :\mathsf{H}^{-\nicefrac12}(\Gz) \to \mathsf{H}^{\nicefrac12}(\Gz)$, associated to $\SL_{\tt G}$, given by $\SLO_{\tt G} \nu=\Dtrace{0}{\SL_{\tt G} \nu}$ for any $\nu \in \mathsf{H}^{-\nicefrac12}(\Gz)$. Note that $\SLO_{\tt G} \nu = \Dtrace{i}{\SL_{\tt G} \nu}$ on $\Gi$ for each $i\in\IntM$ and $\nu \in \mathsf{H}^{-\nicefrac12}(\Gz)$. 

\medskip

It follows from potential theory that the function $w$ satisfying~\eqref{eq:PDE_homogen} can be represented by a global single layer potential:
\begin{equation}
\label{eq:GlobRep}
w = \SL_{\tt G} \nu
\end{equation}
for a unique $\nu \in \mathsf{H}^{-\nicefrac12}(\Gz)$ and it holds that
\[
\nu_i = \nu|_{\Gi} = \frac{1}{4\pi} \left( \nabla w\cdot n_0 + \nabla w\cdot n_i\right)
\qquad\mbox{on each $\Gi$}.
\]
Equation~\eqref{eq:GlobRep} allows the unknown function $w$ to be represented more efficiently in terms of $\nu$. Indeed, the function $\nu$ is only defined on the surface $\Gz$.

\medskip

We also consider the global single layer boundary operator $\SLO_{\tt G} :\mathsf{H}^{-\nicefrac12}(\Gz) \to \mathsf{H}^{\nicefrac12}(\Gz)$, associated to $\SL_{\tt G}$, given by $\SLO_{\tt G} \nu=\Dtrace{0}{\SL_{\tt G} \nu}$. Since $w$ is continuous across $\Gi$, we have
\begin{equation}
\label{eq:GPB}
\lambda_i
=
\Dtrace{i}{w}
=
\Dtrace{0}{w}
=
\Dtrace{0}{\SL_{\tt G} \nu}
=
\SLO_{\tt G} \nu
\qquad \mbox{on $\Gi$}.
\end{equation}

\subsection{Integral equation} \label{sec:int_eq}

We now derive an integral equation which uniquely determines the trace $\lambda = \gamma_0(w)$ of $w$ on $\Gamma_0$.

\medskip

For all $i \in \IntM$, from~\eqref{eq:NeuTrInt} and~\eqref{eq:NeuTrExt}, it holds that
$$
\nabla w\cd\ni = \left(2\pi \SLO_i^{-1} + \frac{\epsilon_i}{2r_i} \right) \lambda_i \quad \mbox{ on $\Gamma_i$},
$$
with
$$
\epsilon_i = -1 \ \ \text{if $i = 1, \ldots, M$}, \qquad \epsilon_\infty = 1.
$$
Besides, using the jump condition on the last line of~\eqref{eq:PDE_homogen}, we obtain that, for each $i\in \IntM$, the outward pointing co-normal derivative of $w$ on $\Gi$ (with respect to $\Oz$) is given by
\[
\grad w\cd \nz
= -p\cd n_0 -\frac{\A_i}{\A_0}(\grad w\cd \ni+p \cd n_i) \quad \mbox{ on $\Gi$}. 
\]
On the other hand, we already noticed in~\eqref{eq:GlobRep} that the solution $w$ can be globally represented by the density $\nu$. We therefore have, on each $\Gi$, that
\begin{align*}
\nu_i
=
\frac{1}{4\pi}( \grad w \cd \ni + \grad w \cd \nz)
=
\frac{\A_0-\A_i}{4\pi \A_0}(\grad w\cd \ni +p \cd n_i)
=
\frac{\A_0-\A_i}{4\pi \A_0}\left(\left(2\pi \SLO_i^{-1} + \frac{\epsilon_i}{2r_i} \right) \lambda_i+p\cd n_i\right).
\end{align*}
In view of~\eqref{eq:GPB}, we derive the globally coupled integral equation 
\begin{align}
\label{eq:GlobProb}
\forall i \in \IntM, \quad
\nu - \frac{\A_0-\A_i}{4\pi \A_0}\left(2\pi\SLO_i^{-1} +\frac{\epsilon_i}{2r_i} \right) (\SLO_{\tt G} \nu)|_{\Gamma_i}
= 
\frac{\A_0-\A_i}{4\pi \A_0} p\cd n_i \quad \mbox{ on $\Gi$}.
\end{align}
Note that the global operator $\SLO_{\tt G}$ couples all the local components $\nu_i = \nu|_{\Gamma_i}$ for $i\in \IntM$.

\medskip

For all $i\in \IntM$, denote by $\mathcal{L}_i: \mathsf{H}^{\nicefrac12}(\Gi) \to \mathsf{H}^{-\nicefrac12}(\Gi)$ the operator defined as follows: for any $\widetilde{\lambda}_i \in \mathsf{H}^{\nicefrac12}(\Gi)$,
\begin{equation}
\label{eq:def_Li}
\mathcal{L}_i\widetilde{\lambda}_i = \frac{\A_0-\A_i}{4\pi \A_0}\left(2\pi\SLO_i^{-1} +\frac{\epsilon_i}{2r_i} \right) \widetilde{\lambda}_i \quad \mbox{ on $\Gamma_i$}.
\end{equation}
We also denote by $\mathcal{L}: \mathsf{H}^{\nicefrac12}(\Gz) \to \mathsf{H}^{-\nicefrac12}(\Gz)$ the operator defined as follows: for any $\widetilde{\lambda} \in \mathsf{H}^{\nicefrac12}(\Gz)$ and all $i\in \IntM$, 
$$
\left( \mathcal L \widetilde{\lambda}\right)|_{\Gi} = {\mathcal L}_i \left({\widetilde \lambda}|_{\Gi} \right).
$$
Lastly, let $g\in \mathsf{H}^{-\nicefrac12}(\Gz)$ be the function defined by 
\begin{equation}\label{eq:tpl}
g|_{\Gamma_i} = \frac{\A_0-\A_i}{4\pi \A_0} p\cd n_i \quad \mbox{ for each $i\in \IntM$}. 
\end{equation}
Equation~\eqref{eq:GlobProb} can then be equivalenty rewritten as
\begin{equation}\label{eq:eqnu}
(Id - \mathcal L\SLO_{\tt G})\nu = g.
\end{equation}
To obtain an equation on the trace $\lambda$ of $w$ on $\Gamma_0$, we apply the operator $\SLO_{\tt G}$ to~\eqref{eq:eqnu}, which leads to
\begin{equation}\label{eq:IntRep}
\left( Id - \SLO_{\tt G} \mathcal L \right)\lambda = \SLO_{\tt G}g.
\end{equation}
This is the integral equation that we discretize in what follows.

\section{Discretization}
\label{sec:disc}

After a short review of the basic properties of real spherical harmonics, we introduce a discretization scheme for~\eqref{eq:IntRep} based on truncated series of real spherical harmonics.
 
\subsection{Real spherical harmonics}

We denote by $(\Ylm)_{\ell \in \N, \, -\ell \le m \le \ell}$ the set of real spherical harmonics (for the unit sphere $\mathbb S^2$ of $\R^3$), normalized in such a way that
$$
\langle \Ylm, \Ylpmp \rangle_{\mathbb S^2}	
=
\int_{\mathbb S^2} \Ylm(s) \, \Ylpmp(s) \, ds
= 
\int_0^\pi \int_{-\pi}^\pi \Ylm(\theta,\phi) \, \Ylpmp(\theta,\phi) \, \sin\theta \, d\theta \, d\phi
=
\delta_{\ell \ell'} \delta_{mm'},
$$
where $\delta_{nm}$ denotes the Kronecker symbol. Spherical harmonics can be defined on the surface $\partial \B_r(x_0)$ of a sphere of center $x_0$ and radius $r$ by
\[
\Ylm\left(\frac{\cdot-x_0}{r}\right).
\]
For any $u,v \in L^2(\partial B_r(x_0))$, we define the scaled inner product
\begin{equation}\label{eq:SIP}
\left\langle u,v \right\rangle_{\partial \B_r(x_0)}
=
\frac{1}{r^2} \int_{\partial B_r(x_0)} u(s)\,v(s) \,ds
=
\int_{\mathbb S^2} u(x_0+r s')\,v(x_0+r s') \,ds'.
\end{equation}
The set of spherical harmonics on $\partial \B_r(x_0)$ forms an orthonormal basis of $L^2(\partial B_r(x_0))$, endowed with this scaled inner product. In particular,
\[
\left\langle \Ylm\left(\frac{\cdot-x_0}{r}\right), \Ylpmp\left(\frac{\cdot-x_0}{r}\right) \right\rangle_{\partial \B_r(x_0)}
=
\langle \Ylm, \Ylpmp \rangle_{\mathbb S^2}	
=
\delta_{\ell\ell'} \delta_{mm'}.
\]
Note that the purpose of the scaled inner product is to avoid to scale the basis functions by the factor $1/r$. This allows us to have the same set of basis functions on all spheres.

\medskip

Standard Sobolev norms can easily be expressed using the decomposition of a function on the spherical harmonics basis. More precisely, for any $u \in L^2(\partial B_r(x_0))$, we have
$$
u(x) = \sum_{\ell=0}^{+\infty} \sum_{m=-\ell}^\ell [\u]_\ell^m \, \Ylm\left( \frac{x-x_0}{r}\right),
$$
with
$$
\forall \ell \in \N, \ \ \forall -\ell \leq m \leq \ell, \quad
[\u]_\ell^m 
=
\left\langle u,\Ylm\left( \frac{\cdot-x_0}{r}\right) \right\rangle_{\partial \B_r(x_0)}
=
\int_{\mathbb S^2} u(x_0+rs) \, \Ylm(s) \, ds.
$$
Besides, for any $u$ and $v$ in $L^2(\partial B_r(x_0))$, we have
\begin{align*}
(u,v)_{\L^2(\partial \B_r(x_0))} 
= 
r^2 \left\langle u , v \right\rangle_{\partial \B_r(x_0)}
=
r^2 \sum_{\ell=0}^{+\infty} \sum_{m=-\ell}^\ell [\u]_\ell^m \, [\v]_\ell^m.
\end{align*}
If the function $u$ belongs to $\H^t(\partial \B_r(x_0))$ for some $t\geq0$, then it holds that
$$
\| u \|_{\H^t(\partial \B_r(x_0))}^2 
= 
r^2 \sum_{\ell=0}^{+\infty}  (1+\ell^2)^t \sum_{m=-\ell}^\ell ([\u]_\ell^m)^2.
$$

\medskip
	
In the sequel, for all $i\in \IntM$, we denote by $\langle \cdot, \cdot \rangle_{\Gi}$ the scaled inner product $\langle \cdot, \cdot \rangle_{\partial B_{r_i}(x_i)}$ and by $\Ylm^i$ the function defined on $\Gi$ by $\dps \Ylm^i(x) = \Ylm\left( \frac{x-x_i}{r_i}\right)$. We also define on $\Gamma_0$ the function $\Ybarlm^i$ by $\Ybarlm^i|_{\Gamma_i} = \Ylm^i$ and $\Ybarlm^i|_{\Gj} = 0$ if $j\neq i$. We can now extend this formalism to functions whose support is included in $\Gamma_0 = \cup_{i \in \IntM} \Gi$. More precisely, for any $u\in L^2(\Gamma_0)$, for all $i\in \IntM$, $\ell \in \N$ and $-\ell \leq m \leq \ell$, we denote by 
$$
[\u_i]_\ell^m := \left\langle u|_{\Gi}, \Ylm\left( \frac{\cdot -x_i}{r_i}\right) \right\rangle_{\Gamma_i}.
$$
We hence have
\[
u|_{\Gi}(x) = \sum_{\ell=0}^{+\infty} \sum_{m=-\ell}^\ell [\u_i]_\ell^m \, \Ylm^i(x) \quad \mbox{ on $\Gi$}
\]
and 
$$
u = \sum_{i\in \IntM} \sum_{\ell=0}^{+\infty} \sum_{m=-\ell}^\ell [\u_i]_\ell^m \ \Ybarlm^i.
$$

\bigskip

After discretization of~\eqref{eq:IntRep}, some numerical quadrature rule is needed to practically compute scalar products of the form~\eqref{eq:SIP}. We therefore introduce the scheme defined by a set $\{s_n,\omega_n\}_{n=1}^{N_g} \subset \left( \mathbb S^2 \times \R_+^\star\right)^{N_g}$ of integration points and weights on the unit sphere, where $N_g \in \mathbb{N}^\star$. We define an approximate scaled inner product as follows: for all $i\in \IntM$ and for any $u,v\in \mathcal C(\Gamma_i)$,
\begin{equation}
\label{eq:DefNumInt}
\left\langle u,v \right\rangle_{\Gi, N_g}
:=  
\sum_{n=1}^{N_g} \omega_n \, u(x_i+r_i s_n) \, v(x_i+r_i s_n).
\end{equation}
In practice, in the numerical tests presented in Section~\ref{sec:Num}, we use the Lebedev quadrature rule~\cite{0953-4075-40-23-004,lebedev1999quadrature} which can integrate products of spherical harmonics exactly up to some degree depending on the number of integration points. We refer to~\cite{Stamm1} for more details.

\subsection{Galerkin approximation}
\label{sec:galerkin}

For any $i\in \IntM$, let $\VNi$ be the set of functions spanned by the spherical harmonics on $\Gi$ of degree $\ell$ lower than or equal to $N$:
\[
\VNi = \text{Span} \left\{ \Ylm^i, \quad 0 \leq \ell \leq N, \ \ -\ell \leq m \leq \ell \right\}.
\]
Any $v \in \VNi$ can thus be written as
$$
v = \sum_{\ell=0}^N \sum_{m=-\ell}^\ell \coeff{v} \, \Ylm^i.
$$
The Galerkin approximation of~\eqref{eq:IntRep} is given by: find $\lambda_N \in \mathsf{H}^{\nicefrac12}(\Gz)$ such that, for each $j\in \IntM$, $\lambda_{N,j} = \lambda_N|_{\Gj}$ belongs to $\VNj$ and
\[
\forall i\in\IntM, \quad \forall v_{N,i}\in \VNi, \qquad
\left\langle \lambda_{N,i} - \SLO_{\tt G} \Lop \lambda_N, v_{N,i} \right\rangle_{\Gi}
= 
\left\langle \SLO_{\tt G}g, v_{N,i} \right\rangle_{\Gi},
\]
with $g$ given by~\eqref{eq:tpl}.

\medskip

In practice, the exact inner product $\langle \cdot, \cdot \rangle_{\Gi}$ has to be replaced by the approximate inner product $\langle \cdot ,\cdot \rangle_{\Gi, N_g}$ for some $N_g\in \N$. This leads to the approximate Galerkin problem: find $\lambda_N \in \mathsf{H}^{\nicefrac12}(\Gz)$ such that, for each $j\in \IntM$, $\lambda_{N,j} = \lambda_N|_{\Gj}$ belongs to $\VNj$ and
\begin{equation}
\label{eq:GalApp}
\forall i\in\IntM, \quad \forall v_{N,i}\in \VNi, \qquad
\left\langle \lambda_{N,i} - \SLO_{\tt G} \Lop \lambda_N, v_{N,i} \right\rangle_{\Gi,N_g}
= 
\left\langle \SLO_{\tt G}g, v_{N,i} \right\rangle_{\Gi,N_g}.
\end{equation}
For each $j\in \IntM$, the function $\lambda_{N,j}$ belongs to $\VNj$. It can thus be represented by a set of coefficients $[\lambda_j]_\ell^m \in \R$, with $0\leq \ell \leq N$ and $-\ell \leq m \leq \ell$, as
$$
\lambda_N|_{\Gj} = \lambda_{N,j}= \sum_{\ell=0}^N \sum_{m=-\ell}^\ell [\lambda_j]_\ell^m \, \Ylm^j.
$$
Therefore the function $\lambda_N$ has a discrete representation in terms of a vector $\lambdav\in \real^{(1+M)(N+1)^2}$ collecting all the coefficients $[\lambda_j]_\ell^m$. 

\medskip

The linear system corresponding to~\eqref{eq:GalApp} is denoted by 
\begin{equation}
\label{eq:LinSys}
\Kv \lambdav = \fv.
\end{equation}
The solution matrix $\Kv$ and the right-hand side vector $\fv$ are given by
\begin{align}
\label{eq:DefL}
[\Kv_{ij}]_{\ell \ell'}^{mm'} &= \left\langle \delta_{ij} \Ybarlpmp^j - \SLO_{\tt G} \Lop \Ybarlpmp^j, \Ybarlm^i \right\rangle_{\Gamma_i, N_g},
\\
\label{eq:DefF}
[\f_i]_\ell^m &= \left\langle \SLO_{\tt G}g, \Ybarlm^i \right\rangle_{\Gi,N_g}.
\end{align}
We now explicitly compute these values.

\medskip

We start by the computation of $\Kv$. We first observe that, for all $j\in \IntM$, the spherical harmonics are eigenfunctions of the local single layer operator $\SLO_j$:
\[
\SLO_j \Ylpmp^j = \frac{4\pi r_j}{2\ell'+1} \Ylpmp^j \qquad \mbox{on $\Gj$}.
\]
In view of the definition~\eqref{eq:def_Li} of $\Lop_j$, we thus obtain that
\begin{align}
\label{eq:DefLApp1}
\Lop_j \Ylpmp^j
=
\frac{\A_0-\A_j}{4\pi \A_0} \left(2\pi\SLO_j^{-1} + \frac{\epsilon_j}{2r_j}\right) \Ylpmp^j
=
\frac{2\ell'+1+\epsilon_j}{2 r_j} \ \frac{\A_0-\A_j}{4\pi \A_0} \ \Ylpmp^j.
\end{align}
Second, the single layer potential generated by the surface charge distribution $\Ylpmp^j$ on $\Gj$ is given by
\begin{align}
\label{eq:DefLApp2_Int}
\forall x \in \Oj, \qquad (\SL_j \Ylpmp^j)(x)
&=
\frac{4\pi r_j}{2\ell'+1} \left(\frac{|x-x_j|}{r_j}\right)^{\ell'}\Ylpmp \left( \frac{x-x_j}{|x-x_j|}\right),
\\
\label{eq:DefLApp2_Ext}
\forall x\in \mathbb R^3\setminus \Oj, \qquad (\SL_j \Ylpmp^j)(x)
&=
\frac{4\pi r_j}{2\ell'+1} \left(\frac{r_j}{|x-x_j|}\right)^{\ell'+1}\Ylpmp \left( \frac{x-x_j}{|x-x_j|}\right).
\end{align}
We now compute $\SL_{\tt G} \Lop\Ybarlpmp^j$, using~\eqref{eq:def_tildeSG} and~\eqref{eq:DefLApp1}:
\begin{equation}\label{eq:jeudi}
\SL_{\tt G} \Lop\Ybarlpmp^j
=
\sum_{i\in\IntM} \SL_i \Lop_i \left( \left. \Ybarlpmp^j \right|_{\Gi} \right)
=
\SL_j \Lop_j \left( \Ylpmp^j \right)
=
\frac{2\ell'+1+\epsilon_j}{2 r_j} \ \frac{\A_0-\A_j}{4\pi \A_0} \ \SL_j \Ylpmp^j.
\end{equation}
In order to compute~\eqref{eq:DefL}, we need to evaluate the function $\SLO_{\tt G} \Lop\Ybarlpmp^j$ (which is the trace of $\SL_{\tt G} \Lop\Ybarlpmp^j$ on $\Gamma_0$) at the points $x_i+r_i s_n$, where we recall that $s_n$ are the integration points on the unit sphere (see~\eqref{eq:DefNumInt}). Introducing
\[
v_n^{ij} := x_i+r_is_n - x_j,
\qquad
s_n^{ij} := \frac{v_n^{ij}}{|v_n^{ij}|}
\qquad \mbox{and} \qquad
t^{ij}_n := \frac{|v_n^{ij}|}{r_j},
\]
and collecting~\eqref{eq:jeudi} with~\eqref{eq:DefLApp2_Int} and~\eqref{eq:DefLApp2_Ext}, we obtain
$$
(\SLO_{\tt G} \Lop \Ybarlpmp^j)(x_i+r_i s_n)
= 
\frac{2\ell'+1+\epsilon_j}{4\ell'+2} \ \frac{\A_0-\A_j}{\A_0} \ (t^{ij}_n)^{f(\ell',\epsilon_j)} \ \Ylpmp(s_n^{ij})
$$
with
\[
f(\ell,\epsilon) = \begin{cases}
\ell & \mbox{if $\epsilon=1$},\\
-(\ell + 1) & \mbox{if $\epsilon=-1$}.
\end{cases}
\]
This leads to
\begin{align*}
\left\langle \SLO_{\tt G} \Lop \Ybarlpmp^j, \Ybarlm^i \right\rangle_{\Gamma_i,N_g}
&=
\sum_{n=1}^{N_g} \omega_n \, \Ylm(s_n) \, (\SLO_{\tt G} \Lop\Ybarlpmp^j)(x_i+r_i s_n)
\\
&= 
\frac{2\ell'+1+\epsilon_j}{4\ell'+2} \ \frac{\A_0-\A_j}{\A_0} \ \sum_{n=1}^{N_g} \omega_n  \, \Ylm(s_n) \, (t^{ij}_n)^{f(\ell',\epsilon_j)} \, \Ylpmp(s_n^{ij}).
\end{align*}
We therefore obtain that
\begin{align}
\label{eq:Lii}
[\Kv_{ii}]_{\ell \ell'}^{mm'} 
&=  
\delta_{\ell \ell'} \delta_{m m'} \left( 1 - \frac{2\ell'+1+\epsilon_i}{4\ell'+2} \ \frac{\A_0-\A_i}{\A_0} \right),
\\
\label{eq:Lij}
[\Kv_{ij}]_{\ell \ell'}^{mm'} 
&= 
- \sum_{n=1}^{N_g} [\Kv_{ij}^n ]_{\ell \ell'}^{m m'} \, (t^{ij}_n)^{f(\ell',\epsilon_j)} \, \Ylpmp(s_n^{ij}), \qquad j\neq i,
\end{align}
with
\[
[\Kv_{ij}^n ]_{\ell \ell'}^{m m'} 
:=
\omega_n \, \frac{2\ell'+1+\epsilon_j}{4\ell'+2} \ \frac{\A_0-\A_j}{\A_0} \ \Ylm(s_n).
\]

\begin{remark}
\label{rem:Lij}
Note that we have analytically computed $[\Kv_{ii}]_{\ell \ell'}^{mm'}$, without resorting to any numerical integration. In contrast, $[\Kv_{ij}]_{\ell \ell'}^{mm'}$ does not have any explicit expression when $j \neq i$, thus the use of a quadrature rule in~\eqref{eq:Lij}. Of course, it is also possible to compute $[\Kv_{ii}]_{\ell \ell'}^{mm'}$ using the numerical integration scheme. We expect to obtain an identical result as soon as a sufficiently large number of Lebedev integration points are used (recall indeed that products of spherical harmonics can be exactly integrated by the Lebedev quadrature rule, provided enough integration points are used). 
\end{remark}

We now turn to the right-hand side $\fv$ of~\eqref{eq:LinSys}. To compute it, we first represent the data appearing in the right hand side of~\eqref{eq:IntRep} by means of spherical harmonics of degree 1: for any $i\in\IntM$, we have
\[
g_i = g|_{\Gi} = \frac{\A_0-\A_i}{4\pi \A_0} p\cd n_i
=
-\epsilon_i \frac{\A_0-\A_i}{4\pi \A_0r_i}(x-x_i)\cd p,
\]
that we recast, for any $x\in\Gi$, as
\begin{align*}
g_i(x)
= 
\sum_{m=-1}^1 [\g_i]_1^m \, \Yom^i(x)
=
\sqrt{\frac{3}{4\pi}} \left[
[\g_i]_1^{-1} \, \frac{(x-x_i) \cdot e_2}{r_i}
+ [\g_i]_1^{0} \, \frac{(x-x_i) \cdot e_3}{r_i}
+ [\g_i]_1^1 \, \frac{(x-x_i) \cdot e_1}{r_i}
\right].
\end{align*}
We thus get
\begin{equation}
\label{eq:ptilde}
[\g_i]_1^m = -\epsilon_i \, \frac{\A_0-\A_i}{4\pi \A_0} \, \sqrt{\frac{4\pi}{3}} \ [\p]_1^m,
\end{equation}
with $[\p]_1^{-1} = p\cdot e_2$, $[\p]_1^{0} = p\cdot e_3$ and $[\p]_1^1 = p\cdot e_1$. We deduce that
$$
\SL_{\tt G}g
=
\sum_{j\in\IntM} \SL_j g_j
=
\sum_{j\in\IntM} \sum_{m'=-1}^1 [\g_j]_1^{m'} \, \SL_j \Yomp^j.
$$
In view of~\eqref{eq:DefF}, we hence have
\begin{align*}
  [\f_i]_\ell^m
  &=
  \sum_{n=1}^{N_g} \omega_n \, (\SLO_{\tt G} g)(x_i+r_i s_n) \, \Ylm(s_n)
  \\
  &=
  \sum_{j\in\IntM} \sum_{m'=-1}^1 [\g_j]_1^{m'} \, \sum_{n=1}^{N_g} \omega_n \, (\SL_j \Yomp^j)(x_i+r_i s_n) \, \Ylm(s_n).
\end{align*}
Using~\eqref{eq:DefLApp2_Int} and~\eqref{eq:DefLApp2_Ext}, we deduce that
\begin{align}
[\f_i]_\ell^m
&=
[\g_i]_1^m \, \frac{4 \pi r_i}{3} \, \delta_{\ell 1}
+
\sum_{\substack{j\in \IntM \\ j\neq i}} \sum_{m'=-1}^1 [\g_j]_1^{m'} \, \sum_{n=1}^{N_g} \omega_n \, (\SL_j \Yomp^j)(x_i+r_is_n) \, \Ylm(s_n)
\nonumber
\\
&=
[\g_i]_1^m \, \frac{ 4\pi r_i}{3} \, \delta_{\ell 1}
+
\sum_{\substack{j\in \IntM \\ j\neq i}} \frac{4\pi r_j}{3} \sum_{n=1}^{N_g} [\h_j^n]_\ell^m \, (t_n^{ij})^{f(1,\epsilon_j)}
\label{eq:vendredi}
\end{align}
with 
\[
[\h_j^n]_\ell^m
:=
\omega_n \, \Ylm(s_n) \sum_{m'=-1}^1 [\g_j]_1^{m'} \, \Yomp\left( s_n^{ij} \right).
\]
We have now explicitly computed the matrix and the right-hand side of the linear system~\eqref{eq:LinSys}. In practice, this linear system is solved with an iterative GMRES solver. In addition, we use the Fast Multipole Method (FMM) to accelerate matrix-vector multiplications when solving~\eqref{eq:LinSys}. Details are postponed until Appendix~\ref{ssec:FMM}. Note that the use of the FMM method is critical in reducing the computational cost.

\subsection{Computation of the energy}

Once~\eqref{eq:LinSys} has been solved for $\lambdav$, we are in position to compute the energy $\cJ^{\mathbb{A}}_{R,p}(\Ainf)$ given by~\eqref{eq:nrj_lambda}. Observing that, on $\Gi$, we have $\jmp{\cA^{\mathbb A,\A_\infty} p \cd n} = (\A_0 \, n_0 + \A_i \, n_i)\cd p = (\A_i - \A_0) n_i \cd p = - 4 \pi \A_0 g_i$, we compute that
\begin{align}
\label{eq:EnergyIntRepDisc}
\cJRNp{p} (\Ainf)
=
\sum_{i=0}^M \frac{\A_i |\Oi| }{|\Bz|}
-
\frac{4\pi \A_0}{|\Bz| }\int_{\Gz} g \, \lambda_N
=
\sum_{i=0}^M \frac{ \A_i |\Oi| }{|\Bz|}
- 
\langle \Psiv, \lambdav \rangle,
\end{align}
where the entries of the vector $\Psiv$ are given by
\begin{equation} 
\label{eq:Psi}
[\Psi_i]_\ell^m = \delta_{\ell 1} \, \frac{4\pi \A_0 r_i^2}{|\Bz|} \, [\g_i]_1^m
\end{equation}
and where, for any $\Psiv$ and $\lambdav$, we have set
\[
\langle \Psiv , \lambdav \rangle 
:=
\sum_{i\in\IntM} \sum_{\ell=0}^N \sum_{m=-\ell}^\ell [\Psi_i]_\ell^m \, [\lambda_i]_\ell^m.
\]
The energy $\cJRNp{p}(\Ainf)$ is the discrete approximation of $\cJ^{\mathbb{A}}_{R,p}(\Ainf)$. In view of~\eqref{eq:somme3}, we finally introduce the discrete approximation of $\cJR^{\mathbb{A}}(\A_\infty)$ as
\begin{align*}
\cJRN(\Ainf) := \frac{1}{3}\sum_{i=1}^3 \cJRNp{e_i}(\Ainf).
\end{align*}

\section{Computation of the discrete homogenized coefficients}\label{sec:coeff}

\subsection{Resolution of the optimization problems}

In Sections~\ref{sec:IE} and~\ref{sec:disc}, we have presented the numerical method we propose to compute the energy~\eqref{eq:nrj_lambda} associated to the embedded corrector problem~\eqref{eq:DiffEq2}, for a fixed value of the exterior diffusion coefficient $\Ainf$. We are now in position to compute an approximation of the homogenized coefficient $\A^\star$, following the arguments of Section~\ref{sec:isotropic}.

Using the discretization method introduced in Section~\ref{sec:disc}, we can define three approximate homogenized coefficients, in the spirit of~\eqref{eq:ExCoeffScal1}--\eqref{eq:ExCoeffScal3}:
\begin{align}
\ANR{1} &= \argmax{\Ainf \in [\alpha, \beta]} \cJRN(\Ainf), 
\label{eq:Discopt1_}
\\
\ANR{2} &= \max_{\Ainf \in [\alpha, \beta]} \cJRN(\Ainf) = \cJRN(\ANR{1}),
\label{eq:Discopt2_}
\\
\ANR{3} \in [\alpha, \beta] \ \ & \text{such that} \ \ \ANR{3} = \cJRN(\ANR{3}),
\label{eq:Discopt3_}
\end{align}
where the subscript $N$ represents the degree of spherical harmonics used in the discretization.

\bigskip

The computation of these three quantities requires the resolution of the optimization problem~\eqref{eq:Discopt1_} and of the fixed-point problem~\eqref{eq:Discopt3_}. We propose the following strategy to solve these problems. First, in the numerical tests we have performed, the mapping 
\[
\Ainf \mapsto \cJRN(\Ainf)
\]
shows good contraction properties (see Figures~\ref{fig:Jr} and~\ref{fig:Jc} below).
We therefore suggest to first solve the fixed-point problem~\eqref{eq:Discopt3_} up to some tolerance to obtain $\ANR{3}$. Second, we optimize the energy functional $\cJRN$
using the solution $\ANR{3}$ to~\eqref{eq:Discopt3_} as initial guess. This leads to the computation of $\ANR{1}$ defined by~\eqref{eq:Discopt1_}. Lastly, $\ANR{2}$ is easily obtained from~\eqref{eq:Discopt2_}.

To solve the optimization problem~\eqref{eq:Discopt1_}, we use the Armijo line search since we can compute the derivative of $\cJRN$ based on the solution of some adjoint linear system (see Section~\ref{sec:derivative} below). 

\subsection{Computation of the first derivative} \label{sec:derivative}

In order to compute the derivative of $\cJRN$ with respect to the scalar quantity $\Ainf$, we first differentiate~\eqref{eq:EnergyIntRepDisc}:
\begin{equation}
\label{eq:Force1}
\frac{d\cJRN }{d{\Ainf}}(\Ainf) 
=
-\left\langle \frac{\partial \Psiv}{\partial {\Ainf}}, \lambdav \right\rangle
- \left\langle\Psiv, \frac{\partial \lambdav}{\partial {\Ainf}} \right\rangle.
\end{equation}
Using~\eqref{eq:ptilde}, we observe that 
\[
\frac{\partial [\g_i]_\ell^m}{\partial {\Ainf}}
= 
\delta_{i\infty} \, \delta_{\ell 1} \, \frac{1}{4\pi \A_0} \, \sqrt{\frac{4\pi}{3}} \, [\p]_1^m.
\]
Using the expression of $[\Psi_i]_\ell^m$ given by~\eqref{eq:Psi}, we obtain
\[
\frac{\partial [\Psi_i]_\ell^m }{\partial {\Ainf}} 
=
\delta_{\ell 1} \, \frac{4\pi \A_0 r_i^2}{|\Bz|} \ \frac{\partial [\g_i]_\ell^m}{\partial {\Ainf}}
=
\delta_{i\infty} \, \delta_{\ell 1} \, \frac{r_i^2}{|\Bz|} \, \sqrt{\frac{4\pi}{3}}\, [\p]_1^m.
\]
To compute $\dps \frac{\partial \lambdav}{\partial {\Ainf}}$, we differentiate the linear system $\Kv\lambdav= \fv$ with respect to $\Ainf$:
\[
\frac{\partial \Kv}{\partial {\Ainf}} \lambdav + \Kv \frac{\partial \lambdav}{\partial {\Ainf}} = \frac{\partial \fv}{\partial {\Ainf}},
\]
hence
\[
\frac{\partial \lambdav}{\partial {\Ainf}}
= \Kv^{-1} \hv
\qquad \mbox{with} \qquad
\hv := \frac{\partial \fv}{\partial {\Ainf}} - \frac{\partial \Kv}{\partial {\Ainf}} \lambdav.
\]
We thus obtain that
\begin{equation}
\label{eq:Force2}
\left\langle\Psiv, \frac{\partial \lambdav}{\partial {\Ainf}} \right\rangle
=
\left\langle\Psiv, \Kv^{-1} \hv \right\rangle
=
\left\langle \sv, \hv \right\rangle,
\end{equation}
where $\sv$ is the solution to the adjoint problem $\Kv^\top \sv =\Psiv$. Collecting~\eqref{eq:Force1} and~\eqref{eq:Force2}, we obtain 
\[
\frac{d\cJRN}{d{\Ainf}}(\Ainf) 
=
-\left\langle \frac{\partial \Psiv}{\partial {\Ainf}}, \lambdav \right\rangle
- \left\langle \sv, \hv \right\rangle.
\]
We have mentioned above that we have used the FMM method when solving the primal system~\eqref{eq:LinSys}. Similarly, we have used this method to again accelerate matrix-vector multiplications when solving the adjoint linear system $\Kv^\top \sv =\Psiv$.

\begin{remark}
The use of an adjoint method to compute~\eqref{eq:Force2} becomes more interesting when $\cJRN$ has to be differentiated with respect to {\em several} parameters.
\end{remark}  

\medskip

We are now left with expanding $\hv$ in spherical harmonics. First, taking a closer look at $[\f_i]_\ell^m $, we infer from~\eqref{eq:vendredi} that
\[
\frac{\partial [\f_i]_\ell^m }{\partial {\Ainf}}
=
\frac{\partial}{\partial {\Ainf}} \left[
[\g_i]_1^m \, \frac{ 4\pi r_i}{3} \, \delta_{\ell 1}
+
\sum_{\substack{j\in \IntM \\ j\neq i}} \frac{4\pi r_j}{3} \sum_{n=1}^{N_g} \omega_n \, \Ylm(s_n) \, (t_n^{ij})^{f(1,\epsilon_j)} \, \sum_{m'=-1}^1 [\g_j]_1^{m'} \, \Yomp\left( s_n^{ij} \right) \right].
\]
Using~\eqref{eq:ptilde}, we obtain
\begin{align*}
\frac{\partial [\f_\infty]_\ell^m }{\partial {\Ainf}} 
&=
\delta_{\ell 1} \, \frac{r_\infty}{3\A_0} \, \sqrt{\frac{4\pi}{3}} \, [\p]_1^m,
\\
\frac{\partial [\f_i]_\ell^m }{\partial {\Ainf}}   
&= 
\frac{r_\infty}{3 \A_0} \, \sqrt{\frac{4\pi}{3}} \, \sum_{n=1}^{N_g} \omega_n \, \Ylm(s_n) \, t_n^{i\infty} \sum_{m'=-1}^1 [\p]_1^{m'} \, \Yomp\left( s_n^{i\infty} \right), \qquad i \neq \infty.
\end{align*}
Second, the coefficients of the matrix $\Kv$ can be derived with respect to $\Ainf$, which yields
\begin{align*}
\frac{\partial [\Kv_{ij}]_{\ell \ell'}^{mm'} }{\partial {\Ainf}} 
&= 
0
\quad 
&& 
\forall i\in \IntM, \ j\neq\infty,
\\
\frac{\partial [\Kv_{ij}]_{\ell \ell'}^{mm'} }{\partial {\Ainf}} 
&= 
\frac{1}{\A_0} \, \frac{\ell'+1}{2\ell'+1} \sum_{n=1}^{N_g} \omega_n \, \Ylm(s_n) \, (t^{ij}_n)^{\ell'} \, \Ylpmp\left( s_n^{ij} \right)
\quad 
&&
\forall i=1,\ldots,M, \ j=\infty,
\\
\frac{\partial [\Kv_{ij}]_{\ell \ell'}^{mm'}}{\partial {\Ainf}} 
&=  
\delta_{\ell \ell'} \, \delta_{m m'} \, \frac{1}{\A_0} \, \frac{\ell'+1}{2\ell'+1}
\quad
&& i = j = \infty.
\end{align*}

\section{Numerical results}
\label{sec:Num}

We present here some numerical tests in the deterministic and stochastic homogenization frameworks presented in Section~\ref{sec:homog} to illustrate the performance of the method we propose.

\medskip

Several parameters influence the accuracy of the approximation of the true homogenized coefficient $\A^\star$:
\begin{itemize}
\item First, the truncation of the material by introducing the ball of radius $R$ introduces a model error. 
\item Second, the solution to~\eqref{eq:IntRep} for fixed $R$ is approximated by a Galerkin scheme based on real spherical harmonics of maximum degree $N$ and using numerical quadrature with $N_g$ points. These approximations create a discretization error. In practice, in what follows, we choose $N_g$ such that the product of two spherical harmonics of degree $N$ is exactly integrated.
\item Third, an iterative solver (with a stopping criterion based on some error tolerance on the residual) is used to solve the linear system~\eqref{eq:LinSys}.
\item Finally, the optimization algorithm involved in the computation of the approximate homogenized coefficients~\eqref{eq:Discopt1_}--\eqref{eq:Discopt3_} also requires an error tolerance.
\end{itemize}
In all computations, unless otherwise stated, the following convergence criteria are used:
\begin{itemize}
\item the iterations of the linear solver are stopped as soon as the relative $l^2$ norm of the residual is smaller than $\eta_{\rm ls}=10^{-7}$; 
\item the iterations of the optimization algorithm are stopped when the absolute value of the difference between two consecutive values of the diffusion parameter $\Ainf$ (which also corresponds to the {\em relative error}, since $a^\star$ is of the order of one in the numerical tests below) is smaller than $\eta_{\rm opt}=10^{-5}$.
\end{itemize}
We recall the reader that, unless otherwise stated, we use the heuristic procedure described in Section~\ref{sec:isotropic} to ensure that (for any fixed value of $R>0$) the inclusions $\Oi$, $1 \leq i \leq M$, do not intersect the sphere $\partial \Bz$.

\subsection{Test case 1}

We first consider a case where spherical inclusions of radius $r_n=0.25$ are periodically arranged on the cubic lattice $\mathbb Z^3$. All the spherical inclusions share the same diffusion coefficient $\A_n=10$, while the diffusion constant of the matrix is fixed to $\A_0=1$. Due to the symmetries of the geometrical setting, the value $\cJRNp{e_i}(\Ainf)$ does not depend on $i$.

\medskip

\begin{figure}[t!]
\centering
\includegraphics[scale=0.4]{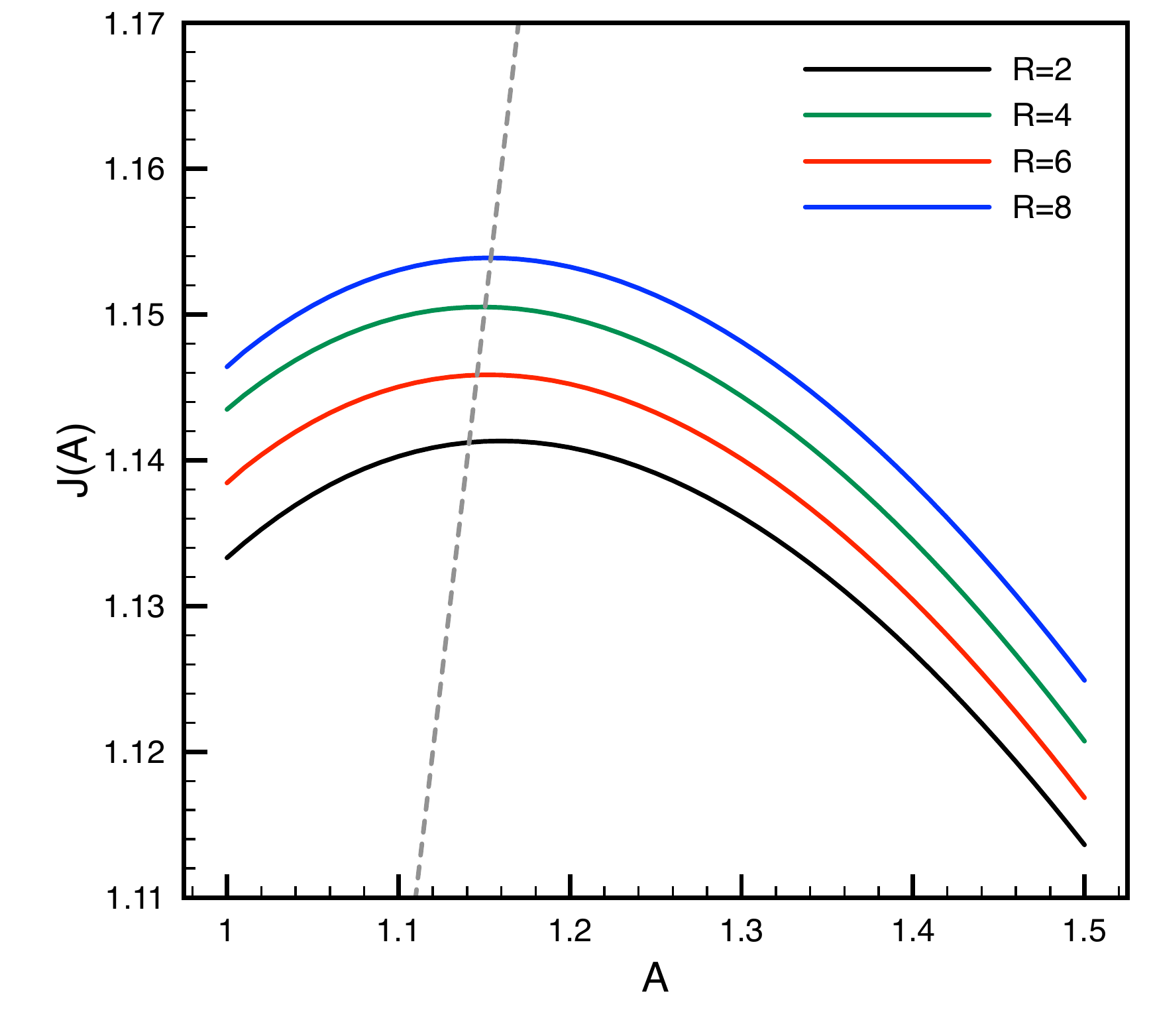}
\caption{[Test case 1] Plot of the function $\Ainf \mapsto \cJRN(\Ainf)$ for $N=1$ and different values of $R$. The dashed line represents the diagonal.}
\label{fig:Jr}
\end{figure}

For the sake of illustration, we show on Figure~\ref{fig:Jr} the function $\Ainf \mapsto \cJRN(\Ainf)$, for different values of $R$ (with $N=1$). As expected, the function is concave. Its derivative is quite small. When $R$ is small (say $R=2$), the value $\ANR{1}$ which maximizes $\cJRN$ is different from the fixed-point value $\ANR{3}$. We expect (if $N$ is sufficiently large) that these two values converge to the same limit when $R \to \infty$. We observe that, when $R$ is large, the values $\ANR{1}$ and $\ANR{3}$ are identical at the scale of the figure.

\medskip

\begin{figure}[t!]
\centering
\includegraphics[scale=0.4]{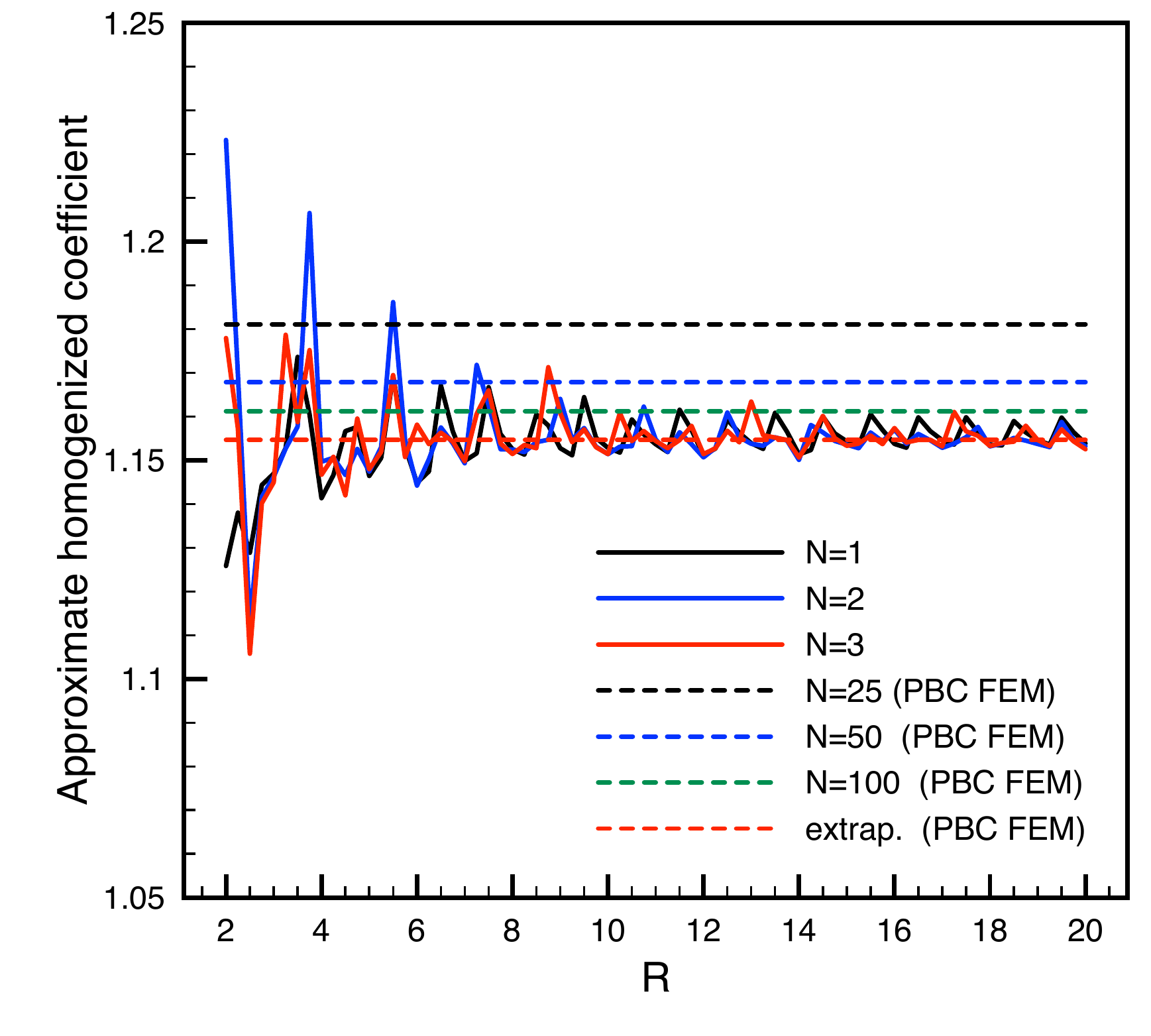}
\includegraphics[scale=0.4]{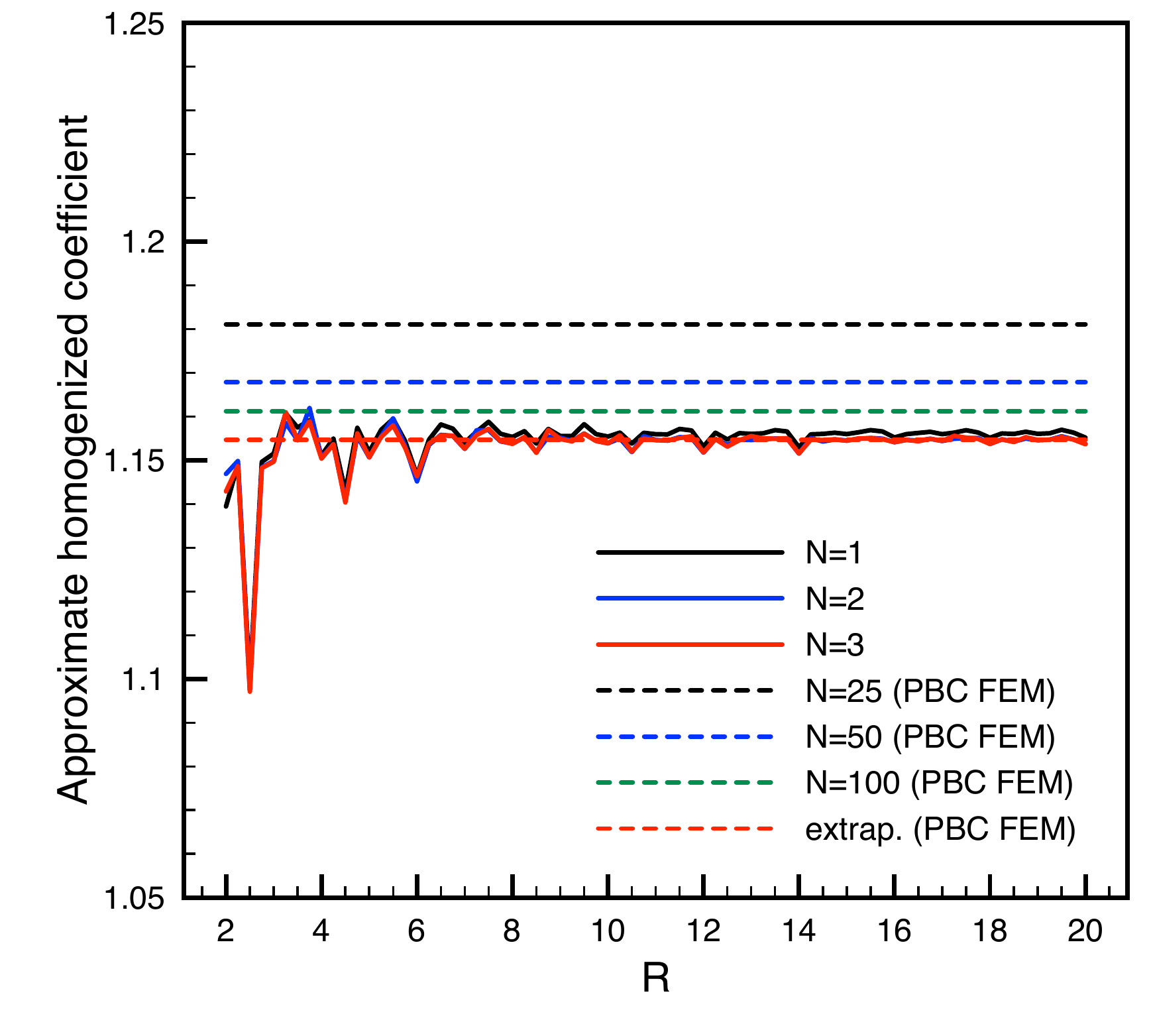}
\caption{[Test case 1] Plots of the functions $R \mapsto \ANR{1}$ (left) and $R \mapsto \ANR{2}$ (right) for different values of the cut-off parameter $N$ (maximal degree of the real spherical harmonics). Dashed lines are the values of the homogenized coefficients obtained from a finite element computation of the ${\mathbb Z}^3$-periodic corrector problem, with $N_d$ mesh points per direction. The coefficients $\ANR{2}$ and $\ANR{3}$ are identical at the scale of the figures.}
\label{fig:T1_long}
\end{figure}

We first investigate the dependency of our approximation with respect to the degree $N$ of the spherical harmonics used in the Galerkin procedure. In Figure~\ref{fig:T1_long}, we show the approximate homogenized coefficient for the values $N=1$, 2, 3 as a function of $R\in[2,20]$. 

We observe that the dependence in $N$ for $\ANR{2}$ and $\ANR{3}$ is negligible with respect to the error in $R$, in the sense that the error due to the truncation in $N$ is dominated by the error introduced by the embedded corrector method on the one hand and by neglecting the spherical inclusions that intersect $\partial \Bz$ on the other hand. The situation is slightly different for $\ANR{1}$ as the dependence in $N$ is more pronounced. The approximations $\ANR{2}$ and $\ANR{3}$ of the exact homogenized coefficient are almost identical, and depend only slightly on $N$. Taking $N=1$ appears to be sufficient.

We have also plotted reference values computed using a finite element code on the unit cell with periodic boundary conditions (PBC). Such comparisons can of course only be made in the simple case of a periodic material. We have plotted the results for three different discretization levels of the finite-element computation (regular mesh with $N_d=25$, $50$ and $100$ discretization points per direction) and have extrapolated the results with the Richardson method. The so-obtained result is regarded as the reference value. Note that such a computation (which is only possible in simple periodic settings) already takes several hours on standard university computer clusters.

\medskip

\begin{figure}[t!]
\centering
\includegraphics[scale=0.4]{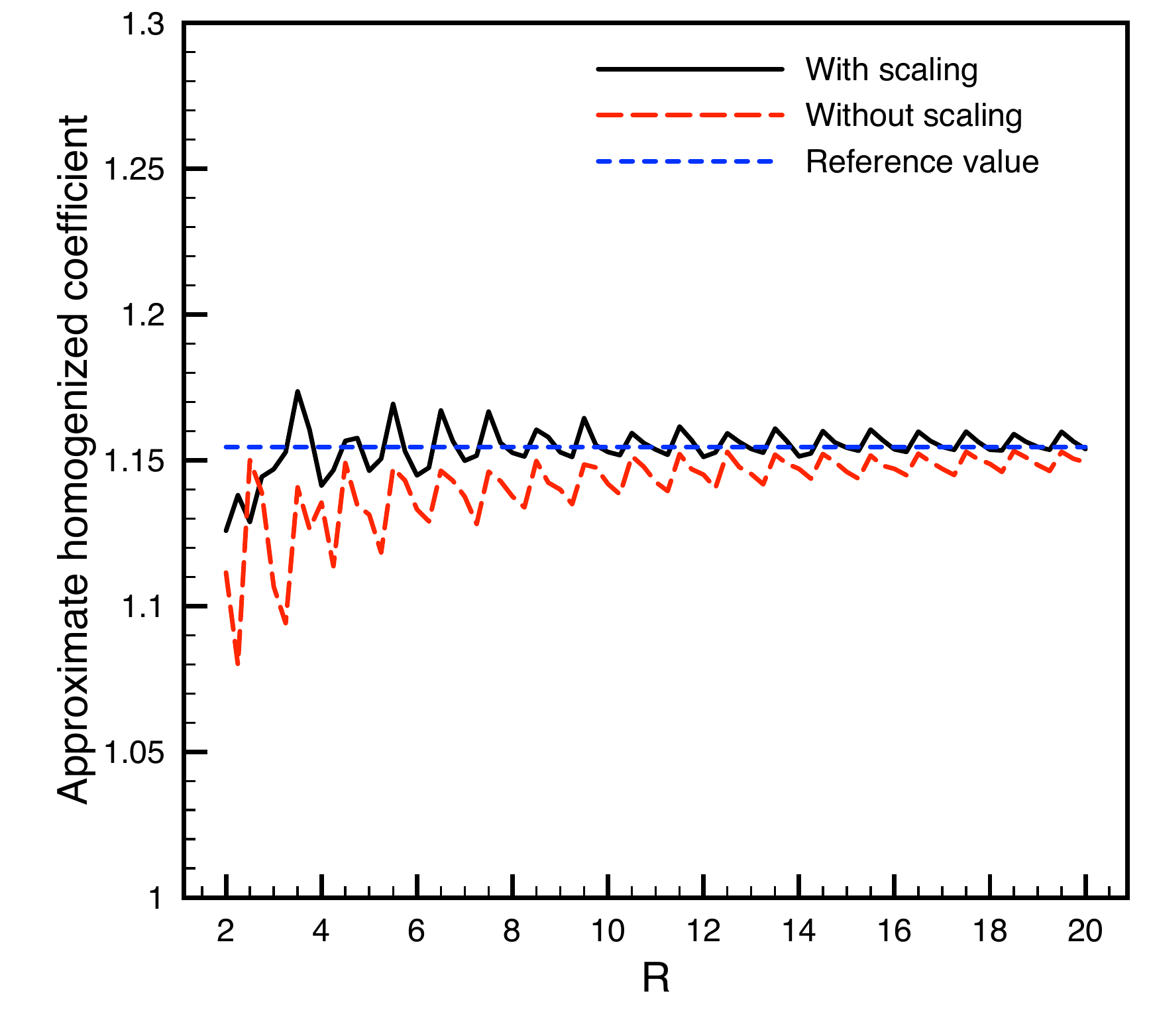}
\includegraphics[scale=0.4]{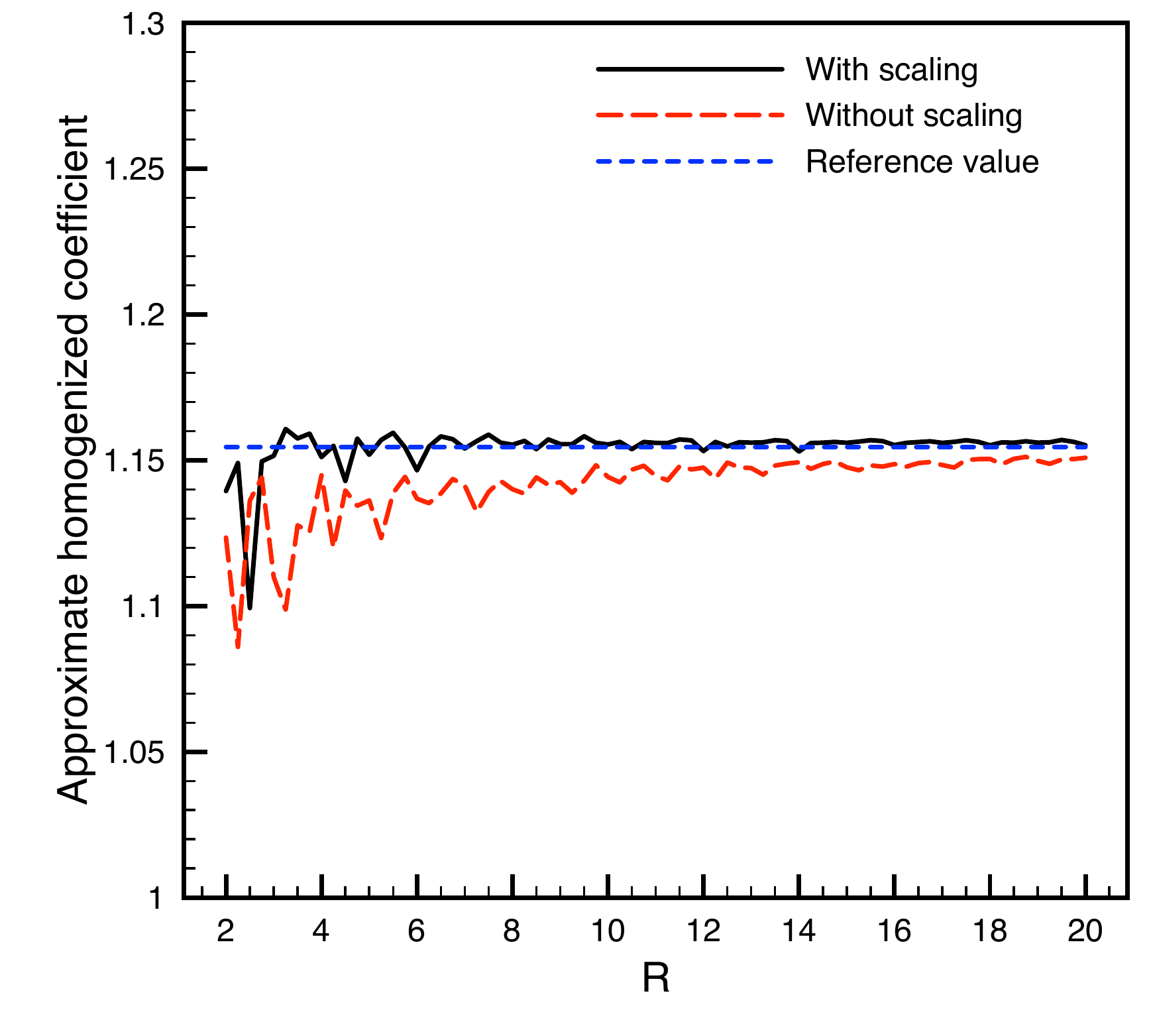}
\caption{[Test case 1] Plots of the functions $R \mapsto \ANR{1}$ (left) and $R \mapsto \ANR{2}$ (right) for $N=1$, with and without scaling of the inclusions inside $\Bz$ (see text). The coefficients $\ANR{2}$ and $\ANR{3}$ are identical at the scale of the figures.}
\label{fig:T1_Scal}
\end{figure}

Figure~\ref{fig:T1_Scal} illustrates how the scaling procedure~\eqref{eq:gamma}--\eqref{eq:scaling}, which is used to better account for the inclusions intersecting $\partial \Bz$, modifies the values of the approximate coefficients $\ANR{1}$, $\ANR{2}$ and $\ANR{3}$. We monitor the homogenized coefficient on the interval $R\in [2,20]$, again with $N=1$.
The dotted line refered to as ``without scaling'' illustrates the value of the effective coefficients computed when the material coefficient inside the ball $\Bz$ is replaced by a coefficient of the form~\eqref{eq:scaling} with $\gamma = 1$ (i.e. when the inclusions intersecting with $\partial \Bz$ have been deleted, but the ones inside $\Bz$ have {\em not} been enlarged). We also report the extrapolated value obtained from the FEM computations as reference value. We observe that the scaling procedure prevents a systematic and very slowly convergent bias of the approximation introduced by discarding the inclusions intersecting $\partial \Bz$. This motivates the scaling procedure~\eqref{eq:gamma}--\eqref{eq:scaling}, which we have used in Figure~\ref{fig:T1_long}, and that we again use in all the following computations, except otherwise stated.

\medskip

\begin{figure}[t!]
\centering
\includegraphics[scale=0.4]{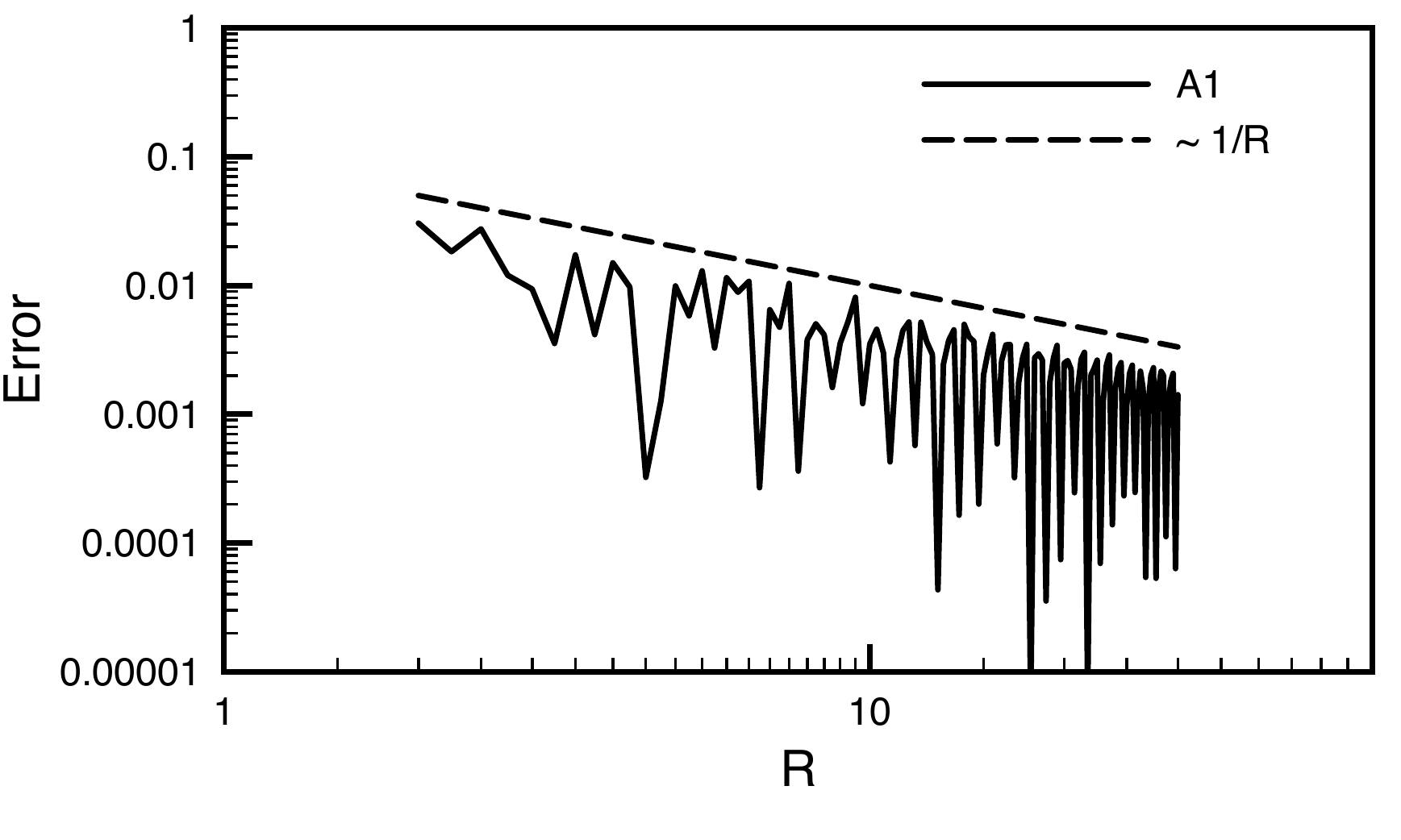}
\includegraphics[scale=0.4]{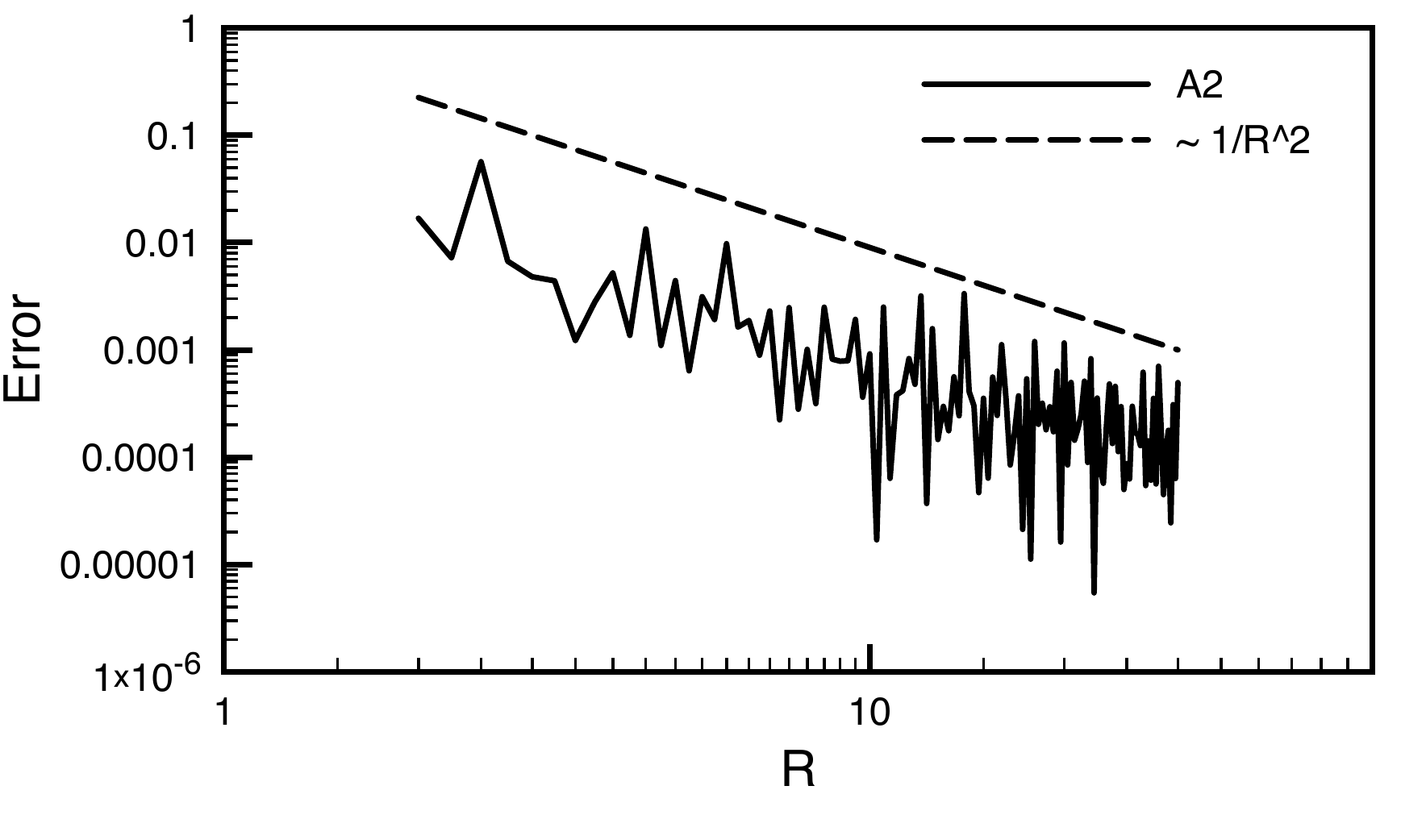}
\caption{[Test case 1] Errors $|\ANR{1}-a^\star|$ (left) and $|\ANR{2}-a^\star|$ (right) as functions of $R$ for $N=1$ (log-log scale). The reference value $a^\star$ has been obtained as explained in the text.}
\label{fig:T1_error}
\end{figure}

To study the convergence with respect to $R$, we have computed a reference solution which is obtained as the average of the results for $R=40$, $40.25$, $40.5$ and $40.75$ with the tighter convergence criteria $\eta_{\rm ls}=10^{-8}$ and $\eta_{\rm opt}=10^{-6}$. This geometrical set-up contains up to 278,370 spherical inclusions (for $R=40.75$). The errors on the coefficients $\ANR{1}$ and $\ANR{2}$ are shown on Figure~\ref{fig:T1_error}. We can observe decays that are proportional at least to $1/R$ for $\ANR{1}$, and approximately to $1/R^2$ for $\ANR{2}$.

\medskip

\begin{figure}[t!]
\centering
\includegraphics[scale=0.4]{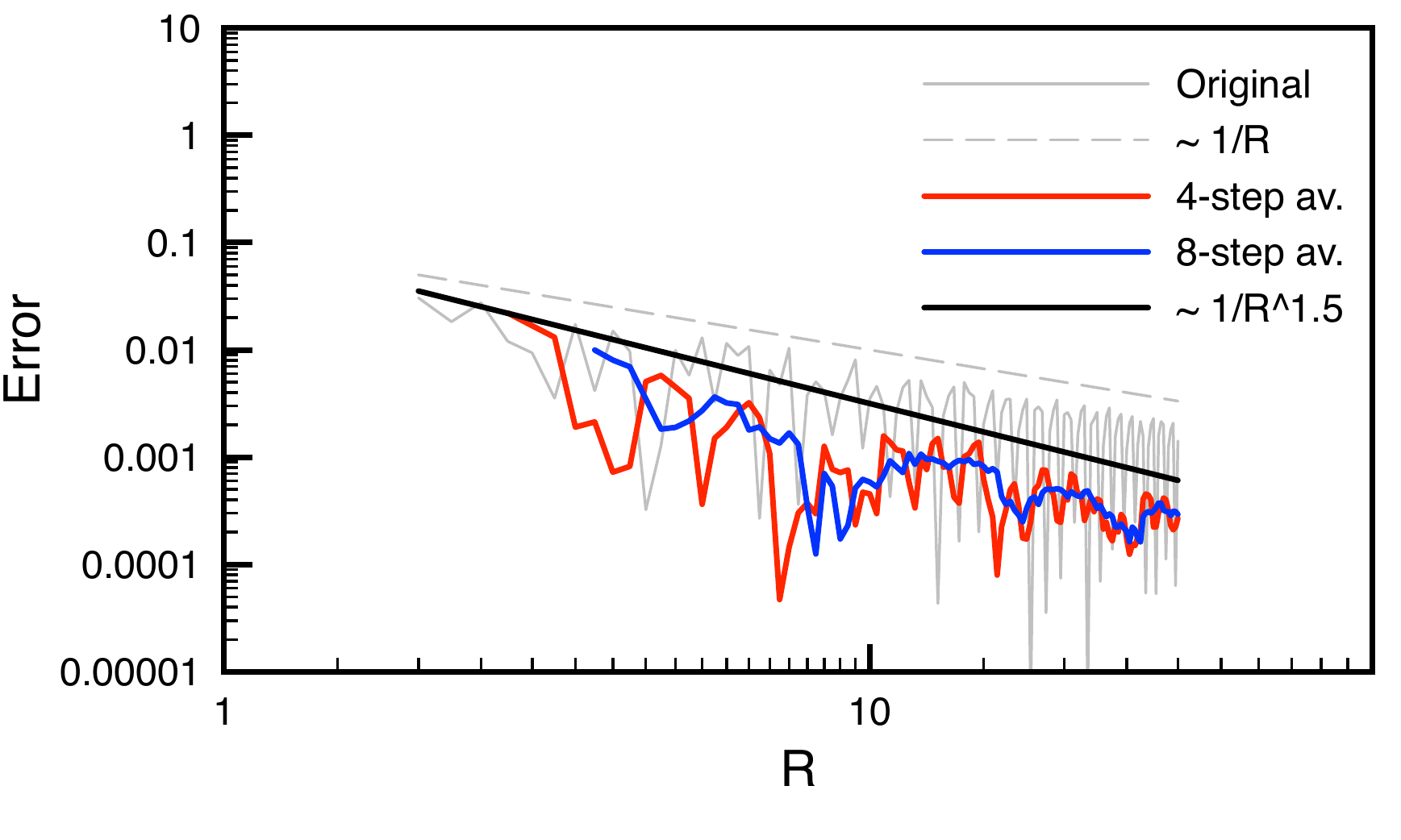}
\includegraphics[scale=0.4]{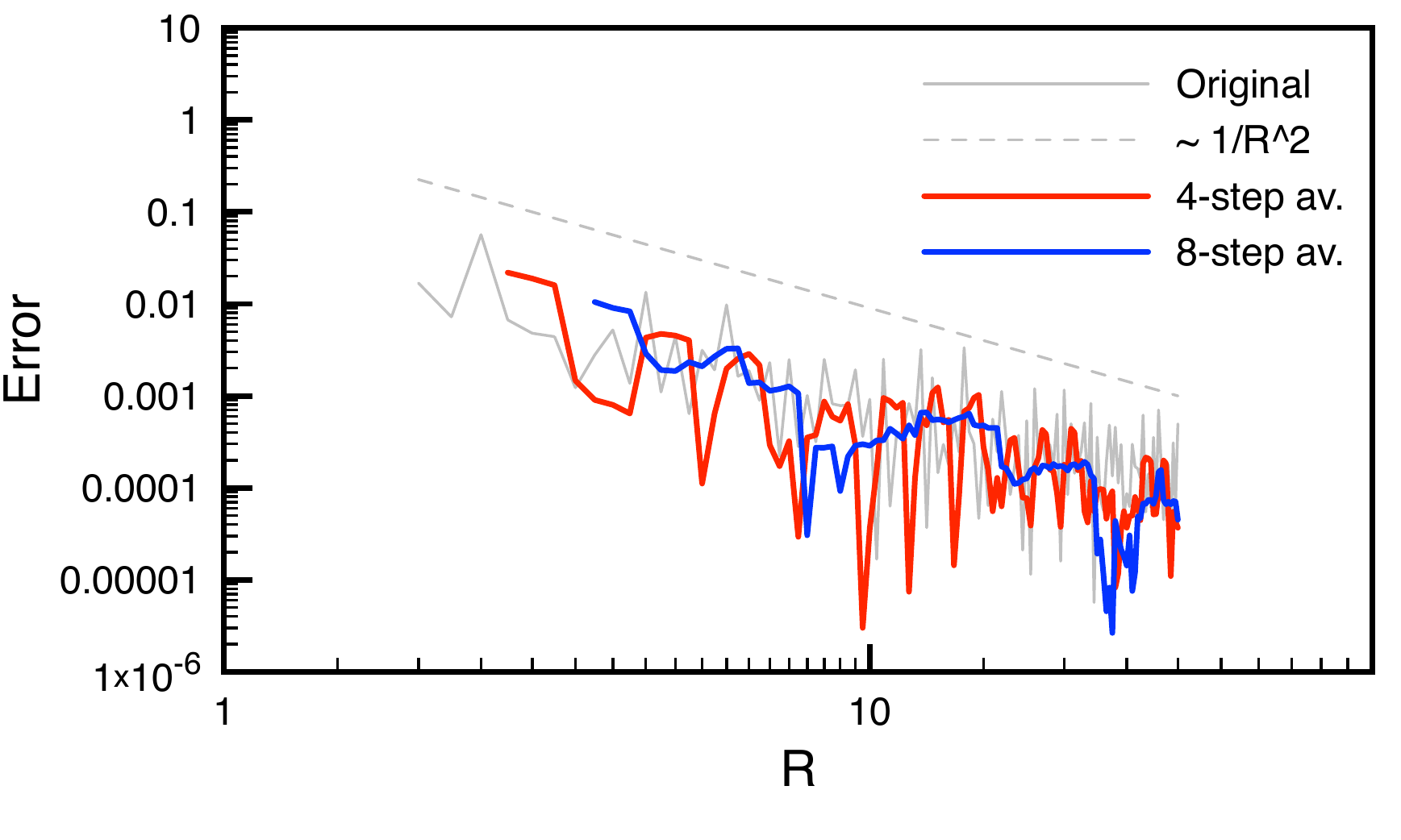}
\caption{[Test case 1] Errors on the approximate homogenized coefficients $\ANR{1}$ (left) and $\ANR{2}$ (right) when local averages with respect to $R$ are performed ($N=1$).}
\label{fig:T1_error_pp}
\end{figure}

The oscillatory nature of the convergence plot and the fact that the local maxima and minima seem to converge according to a different rate indicate that taking averages over different values of $R$ (namely filtering w.r.t.~$R$) might result in a faster convergence. In Figure~\ref{fig:T1_error_pp}, we present convergence plots of approximate homogenized coefficients obtained by averaging over 4 (resp. 8) consecutive values of $R$. The sample points for $R$ are evenly spaced every $0.25$ from $R=2$ to $R=30$. We observe that, in the case of $\ANR{1}$, a better convergence rate (of approximately $R^{-3/2}$ or better) is achieved. The improvement of the convergence rate, if any, is less significant for $\ANR{2}$. 

\subsection{Test case 2}

In this test case, we take the radii of the spherical inclusions smaller, setting $r_n=0.15$, and we increase the contrast between the material coefficients by setting $\A_n=50$. The background coefficient is again set to $\A_0=1$ and the spheres are placed on the cubic lattice $\mathbb Z^3$. All other parameters are identical to the ones of Test case 1.

\medskip

\begin{figure}[t!]
\centering
\includegraphics[scale=0.4]{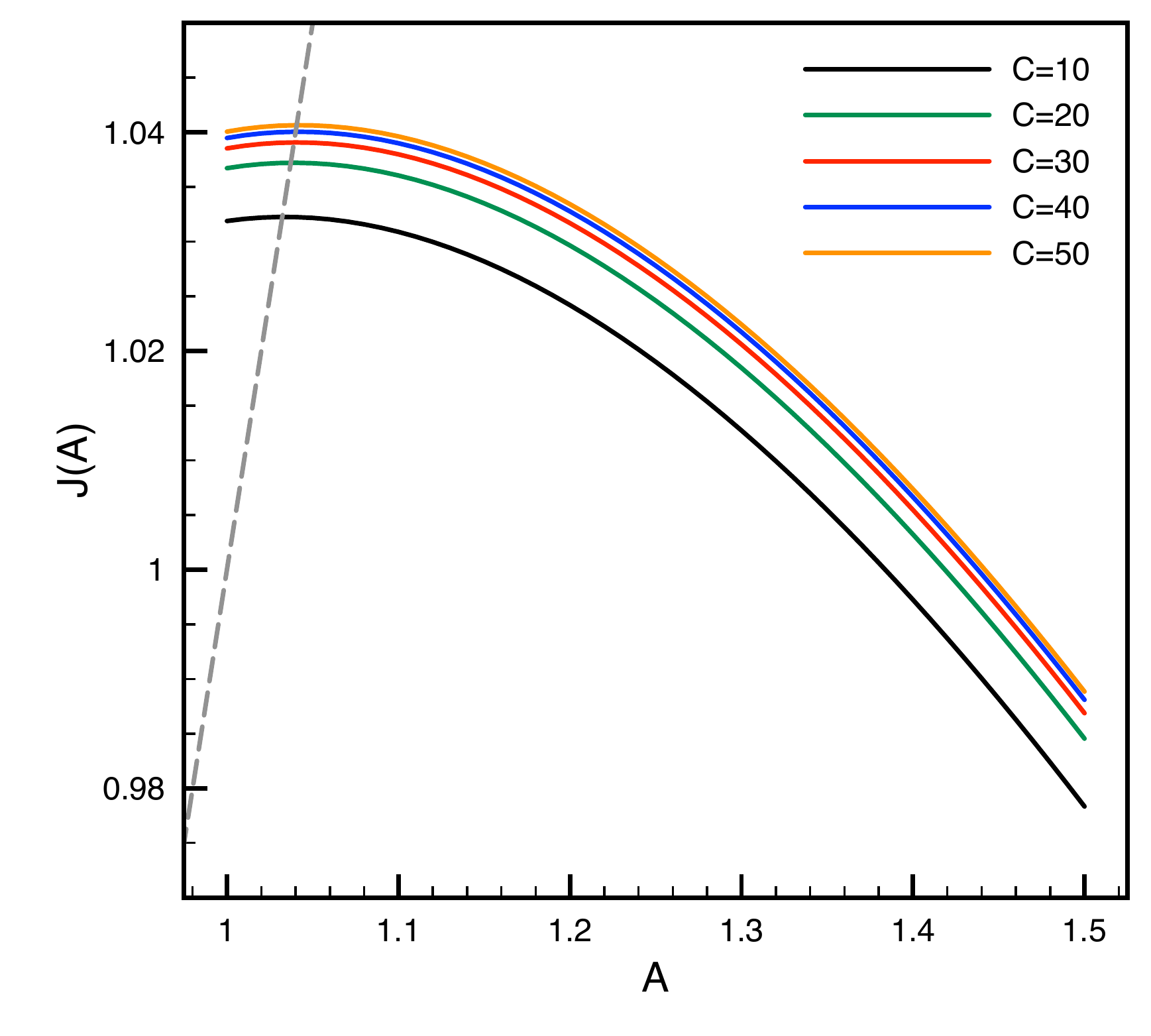}
\caption{[Test case 2] Plot of the function $\Ainf \mapsto \cJRN(\Ainf)$ for $N=1$, $R=6$ and different values of the contrast. The dashed line represents the diagonal.}
\label{fig:Jc}
\end{figure}

For the sake of illustration, we show on Figure~\ref{fig:Jc} the function $\Ainf \mapsto \cJRN(\Ainf)$, for different values of the contrast (with $N=1$ and $R=6$). We thus fix $\A_0=1$ and let $\A_n$ take different values between 10 and 50. We observe that the function depends only mildly upon the contrast. In the sequel of that test, we set $\A_n=50$ as mentioned above.

\medskip

\begin{figure}[t!]
\centering
\includegraphics[scale=0.4]{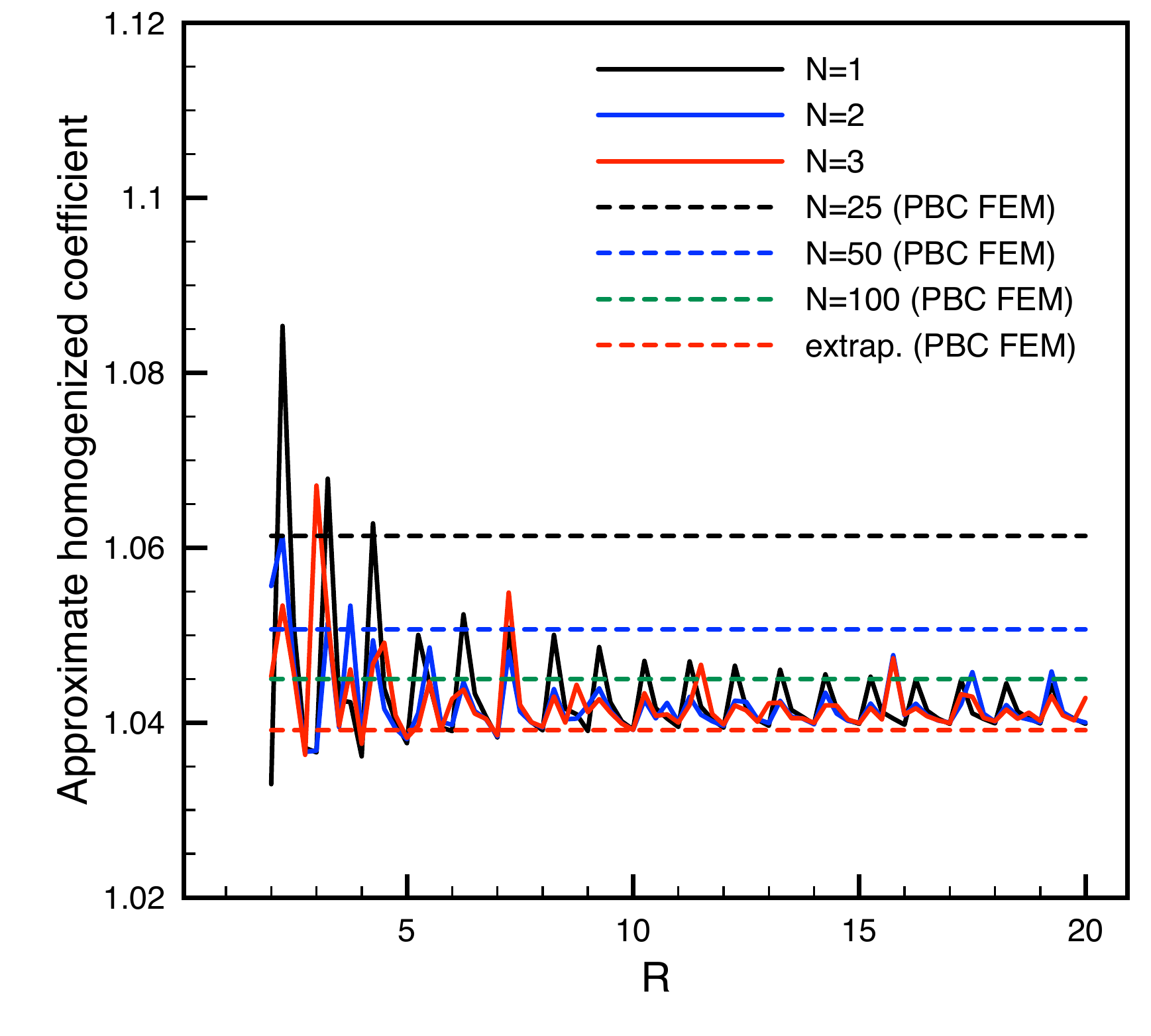}
\includegraphics[scale=0.4]{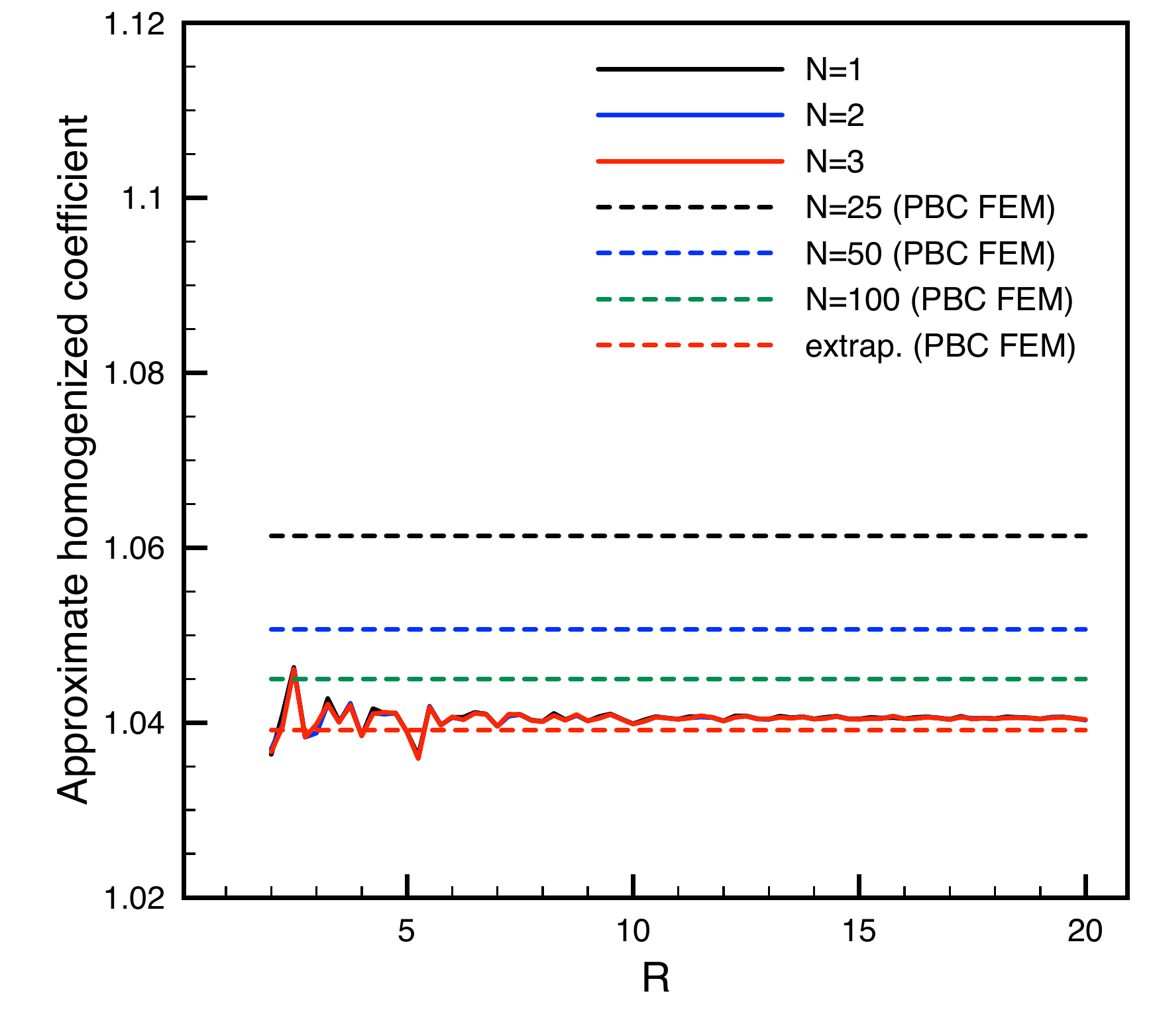}
\caption{[Test case 2] Plots of the functions $R \mapsto \ANR{1}$ (left) and $R \mapsto \ANR{2}$ (right) for different values of the cut-off parameter $N$ (maximal degree of the real spherical harmonics). Dashed lines are the values of the homogenized coefficients obtained from a finite element computation of the ${\mathbb Z}^3$-periodic corrector problem, with $N_d$ mesh points per direction. The coefficients $\ANR{2}$ and $\ANR{3}$ are identical at the scale of the figures.}
\label{fig:T2_long}
\end{figure}

In Figure~\ref{fig:T2_long}, we plot the homogenized coefficients as $R$ varies between 2 and 20. Here, we observe a larger mismatch between the value obtained by a FEM computation (with periodic boundary conditions (PBC) on the unit cell) and our approximation. This might be explained by the fact that this problem is harder to solve by finite element methods. While the method introduced in this article seems to converge for $R>10$ (with a relative error smaller than $10^{-2}$), the FEM does not reach a similar accuracy for the discretization levels we have considered (namely regular meshes with $25$, $50$ and $100$ discretization points per direction).

\begin{figure}[t!]
\centering
\includegraphics[scale=0.4]{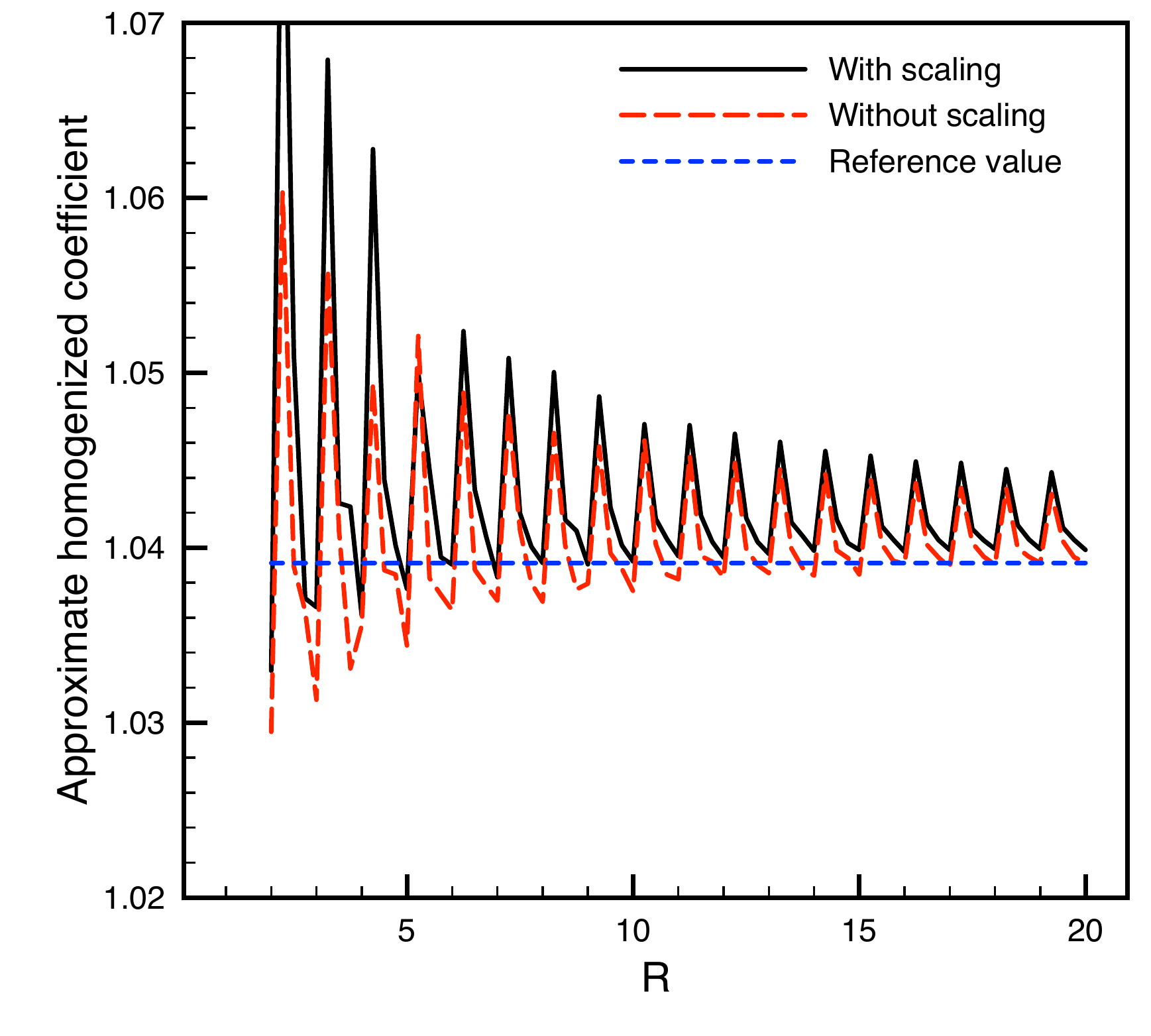}
\includegraphics[scale=0.4]{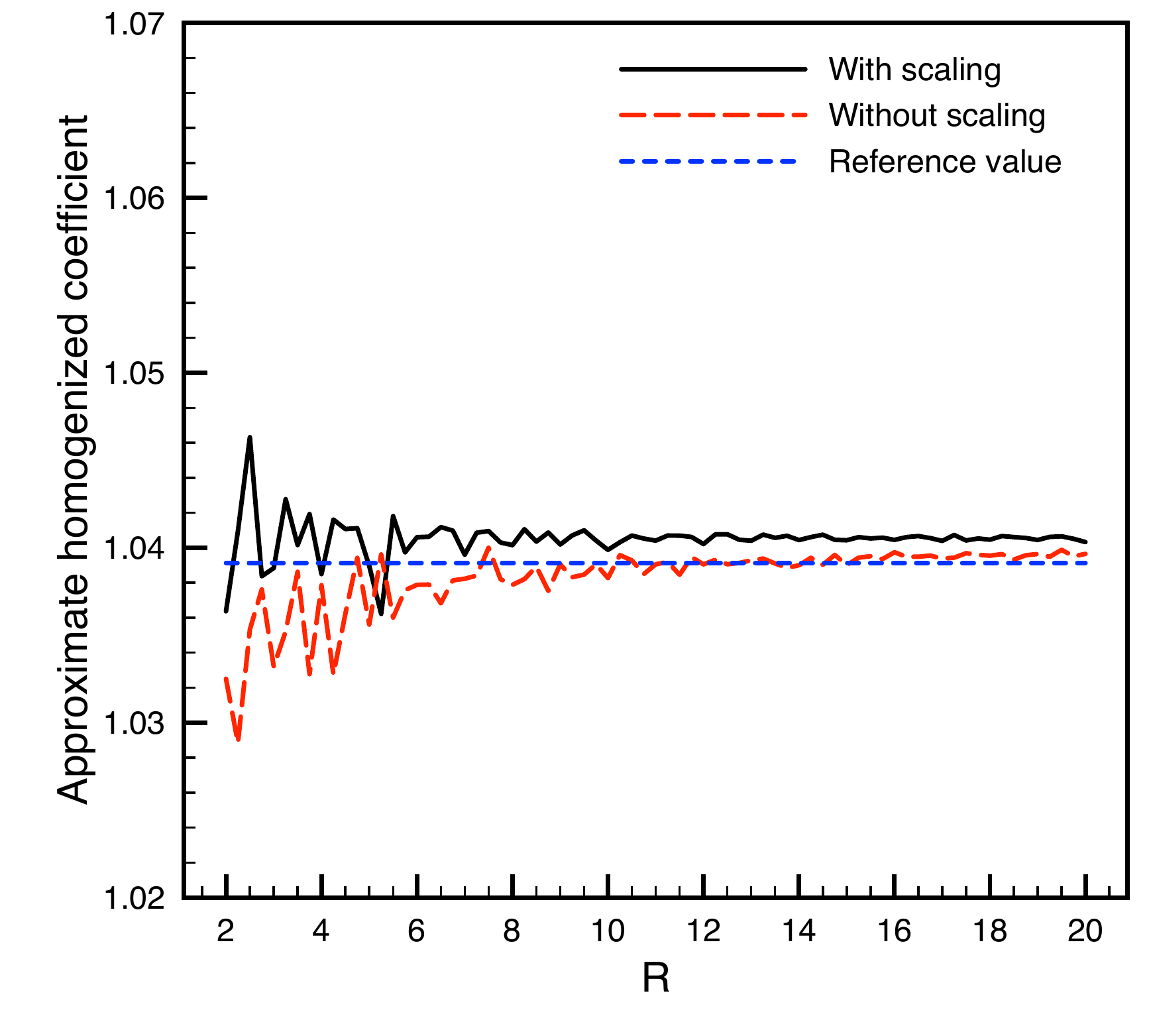}
\caption{[Test case 2] Plots of the functions $R \mapsto \ANR{1}$ (left) and $R \mapsto \ANR{2}$ (right) for $N=1$, with and without scaling of the inclusions inside $\Bz$ (see text). The coefficients $\ANR{2}$ and $\ANR{3}$ are identical at the scale of the figures.}
\label{fig:T2_Scal}
\end{figure}

Figure~\ref{fig:T2_Scal} shows the approximate homogenized coefficients as $R$ varies between 2 and 20, with and without the scaling procedure of the inclusions contained in the ball $\Bz$. In comparison to Test case 1, we observe that the effect of the scaling is weakened for $\ANR{1}$, while it is still clearly visible for $\ANR{2}$ and $\ANR{3}$.

\begin{figure}[t!]
\centering
\includegraphics[scale=0.4]{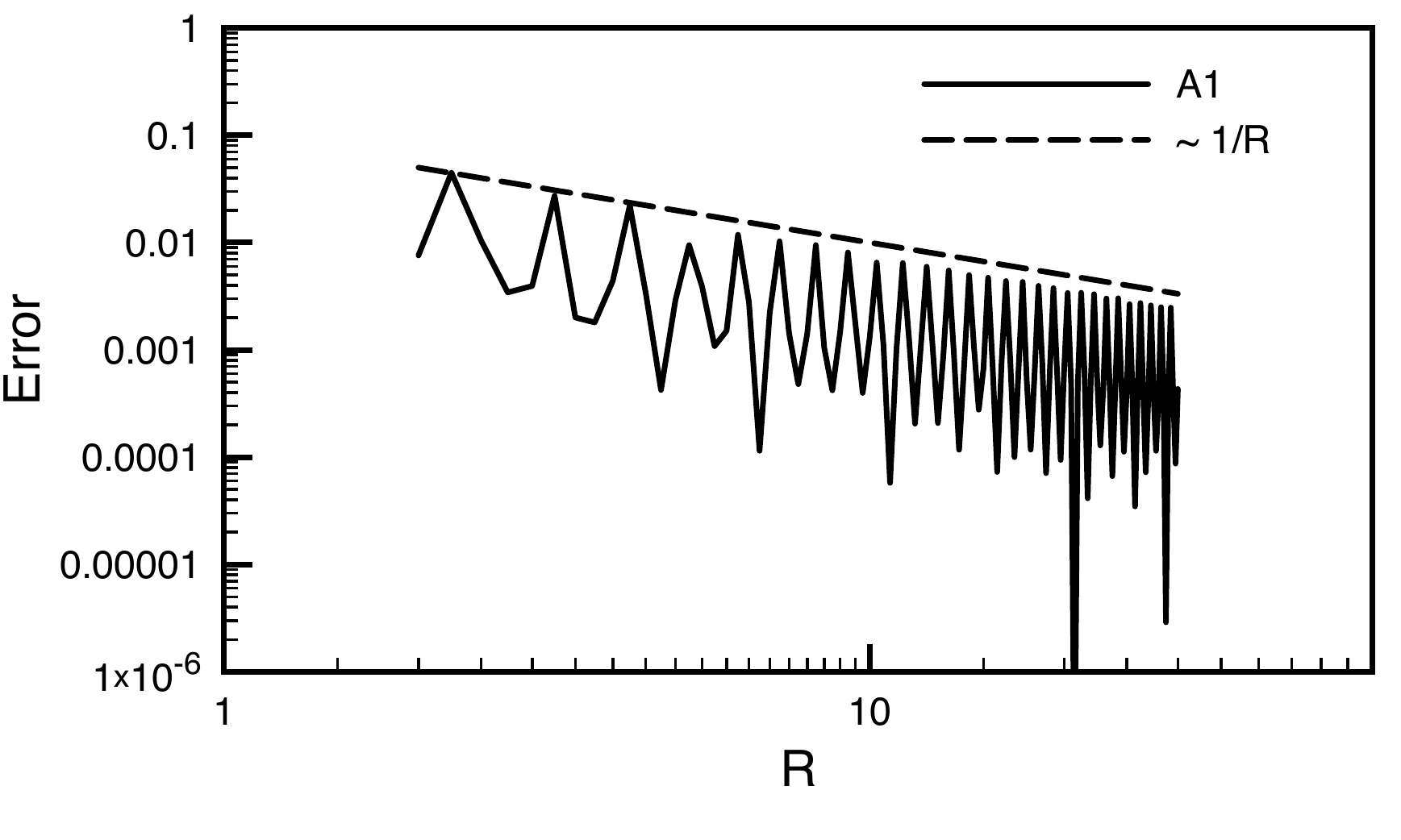}
\includegraphics[scale=0.4]{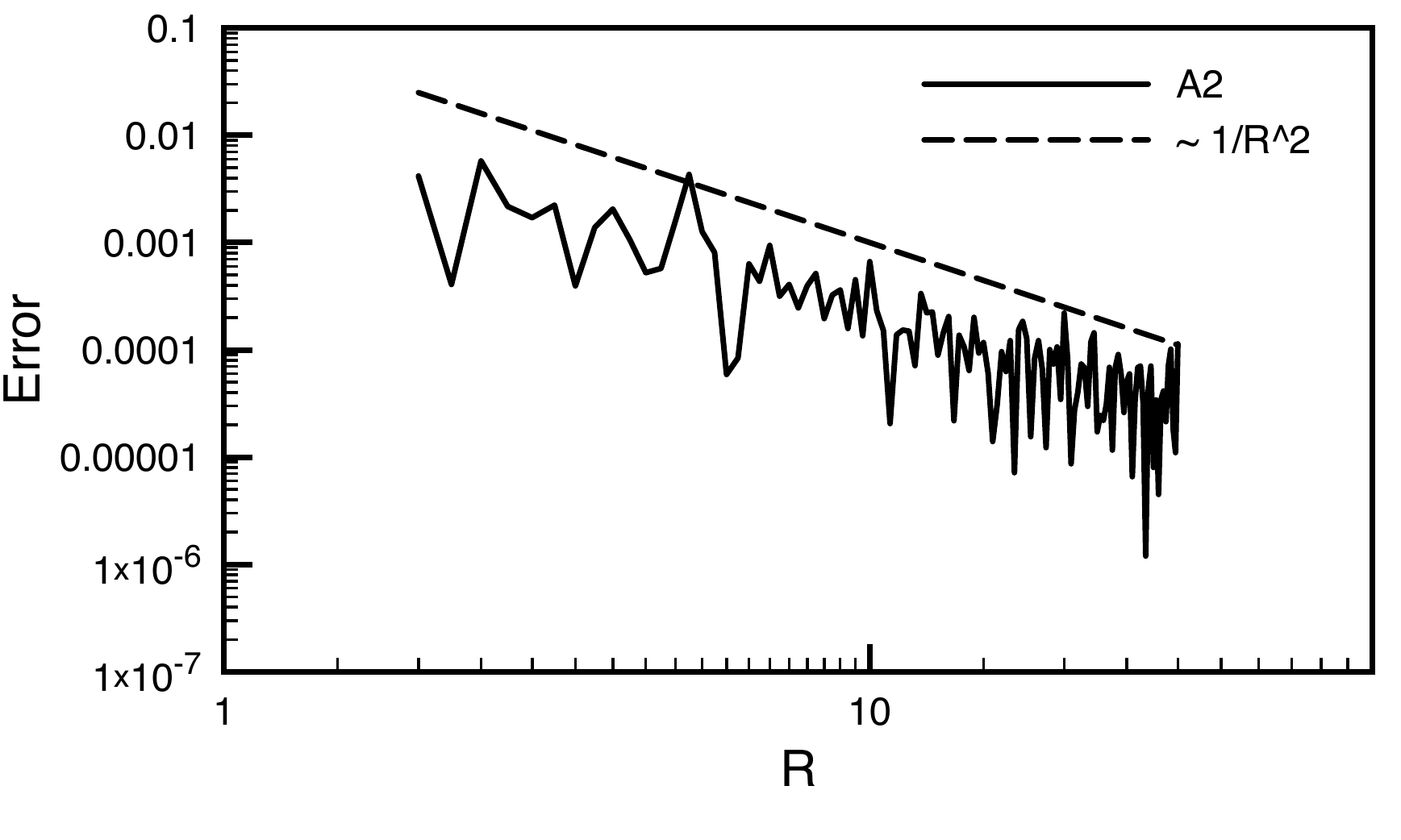}
\caption{[Test case 2] Errors $|\ANR{1}-a^\star|$ (left) and $|\ANR{2}-a^\star|$ (right) as functions of $R$ for $N=1$ (log-log scale). The reference value $a^\star$ has been obtained with the same method as for Figure~\ref{fig:T1_error}.}
\label{fig:T2_error}
\end{figure}

Finally, on Figure~\ref{fig:T2_error}, we study the error on the homogenized coefficients as a function of $R\in[2,30]$. The parameters of this test are the same as the convergence test for Test case 1 (see Figure~\ref{fig:T1_error}). Again, we can observe a convergence rate at least of $1/R$ for $\ANR{1}$ and approximately of $1/R^2$ for $\ANR{2}$.
 
\subsection{Test case 3}

We now consider the case of a random material with polydisperse spherical inclusions. The radii of the inclusions are uniformly distributed between 0.1 and 0.25, their centers are uniformly distributed under the constraint that the distance between two spheres is not smaller than 0.4, and the diffusion coefficient in each inclusion is uniformly distributed between 10 and 50. On average, there is one inclusion per cube of unit size. The three random variables (radius, position and diffusion coefficient) are independent. We consider $R$ in the range $[2,20]$. The largest simulated configuration consists of 32,442 inclusions.

For each computation (with a fixed radius $R$), the property that $\cJRNp{e_i}(\Ainf)$ does not depend on $i$ is of course lost, due to the random positions of the different inclusions, their random radii and their random diffusion coefficients. In the limit $R\to\infty$, we expect however that all $\cJRNp{e_i}(\Ainf)$ converge to the same value (independent of $i$) since the system is statistically isotropic. For simplicity, the numerical results presented in this section have been obtained with the energy functional 
\[
\cJRNp{e_1}(\Ainf)
\qquad \mbox{instead of} \qquad
\frac{1}{3}\sum_{i=1}^3 \cJRNp{e_i}(\Ainf).
\]
Note that the two functionals coincide in the limit $R \to \infty$.

\medskip

\begin{figure}[t!]
\centering
\includegraphics[scale=0.4]{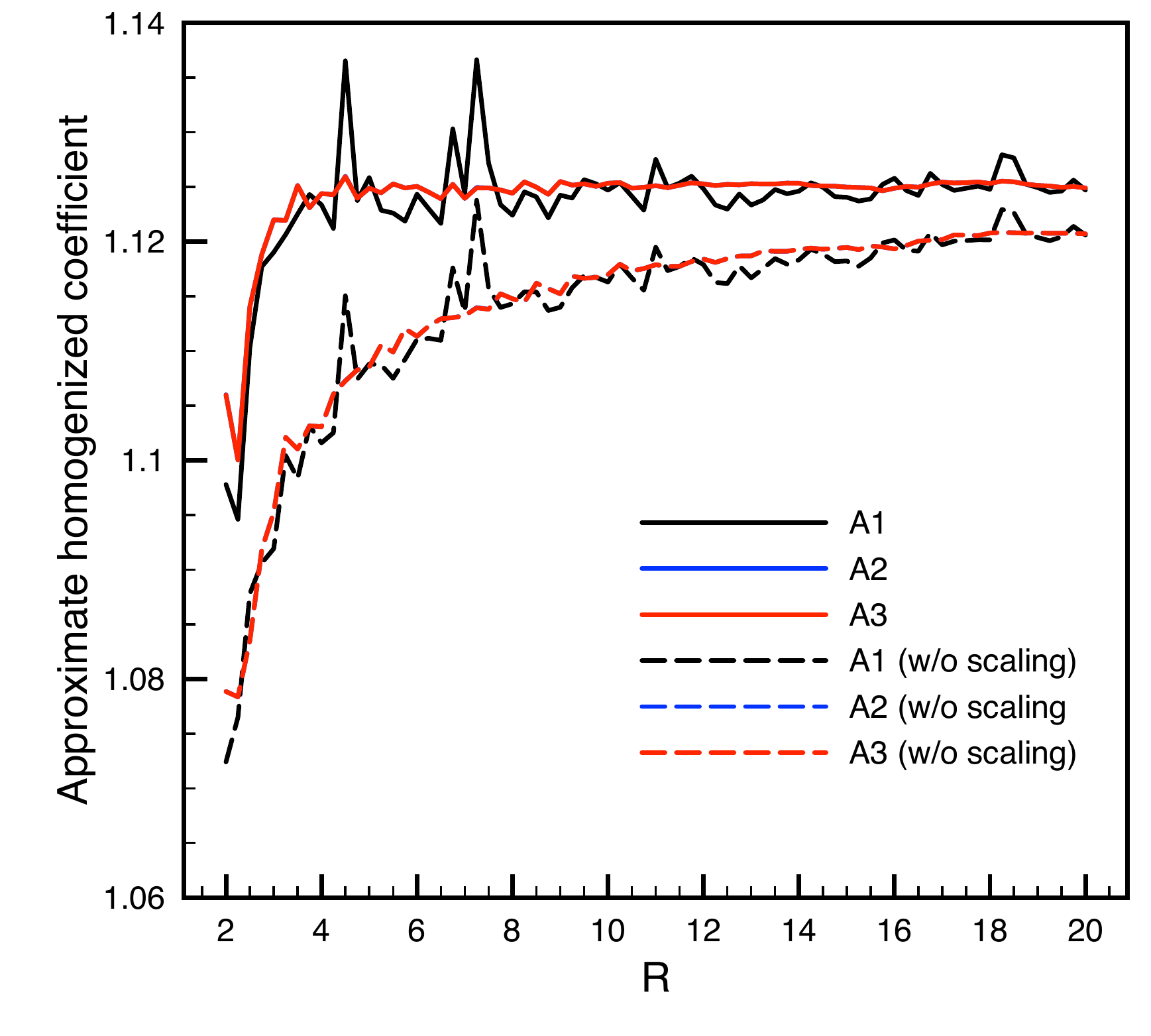}
\includegraphics[scale=0.4]{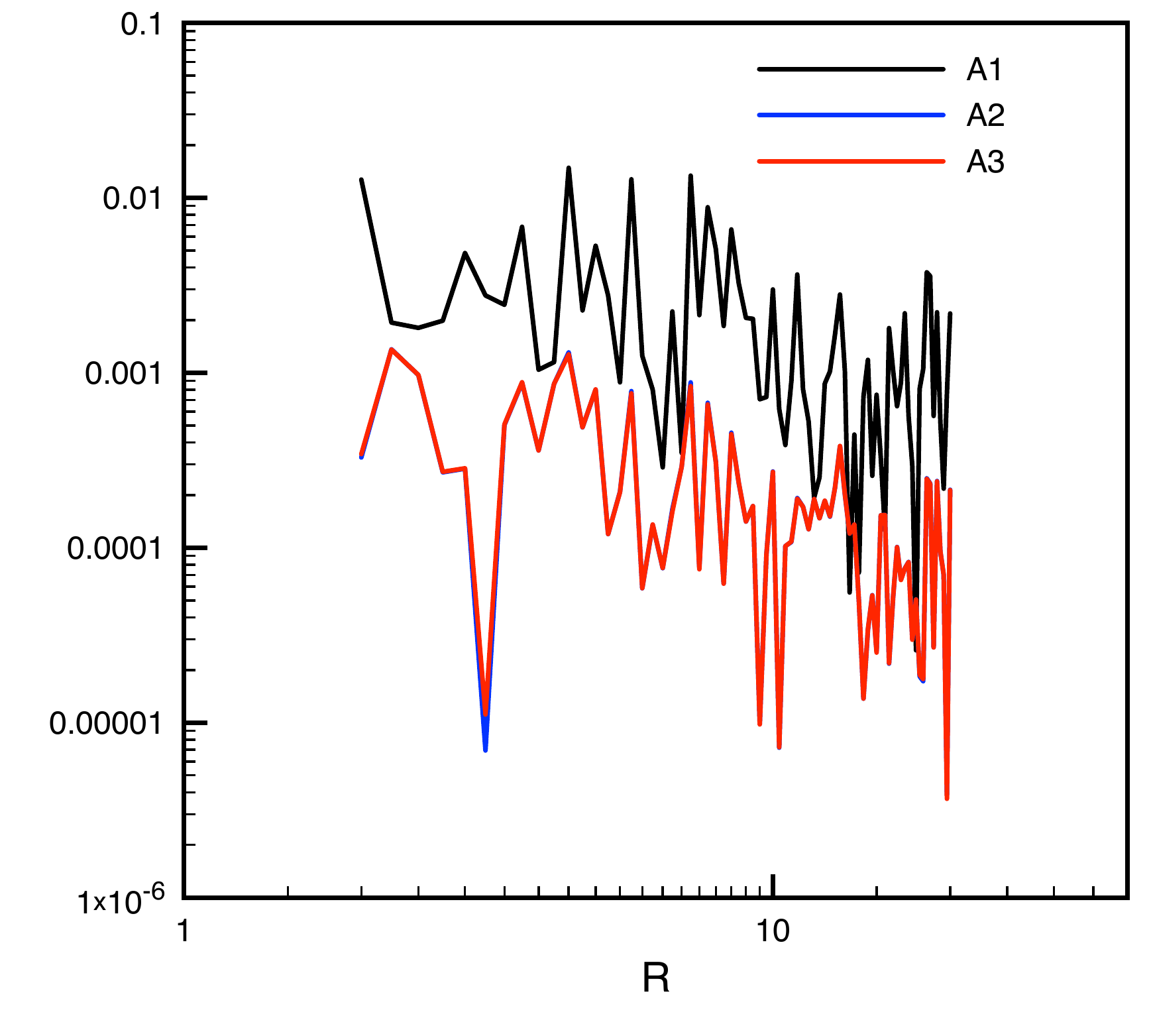}
\caption{[Test case 3] Left: Plots of the functions $R \mapsto \ANR{1}$ (marked as A1), $R \mapsto \ANR{2}$ (marked as A2) and $R \mapsto \ANR{3}$ (marked as A3) for $N=1$. Right: Absolute value of the difference of the homogenized coefficients, when spatial isotropy is assumed and when it is not.}
\label{fig:LongT3}
\end{figure}

On the left side of Figure~\ref{fig:LongT3}, we plot the three approximate homogenized coefficients $\ANR{1}$, $\ANR{2}$ and $\ANR{3}$ as functions of $R$ (we have set $N=1$). We have run two simulations, which are based on the same geometric configuration, that is on the {\em same realization} of the random material. In the first one, we have simply discarded the inclusions intersecting with the boundary $\partial \Bz$, while in the second one the scaling procedure~\eqref{eq:gamma}--\eqref{eq:scaling} is used. The results with and without using the scaling procedure are significantly different, and the approximations converge much faster with respect to $R$ when the scaling procedure is used.

On the right side of Figure~\ref{fig:LongT3}, we show the absolute value of the difference between the two homogenized limits computed when using the functional $\cJRNp{e_1}(\Ainf)$, respectively the functional $\dps \frac{1}{3}\sum_{i=1}^3 \cJRNp{e_i}(\Ainf)$. Despite the fact that, for a finite $R$, the two functionals differ, we observe that the error on the approximate homogenized coefficient is small. As expected, this error decreases when $R$ increases. 

\medskip

\begin{figure}[t!]
\centering
\includegraphics[scale=0.4]{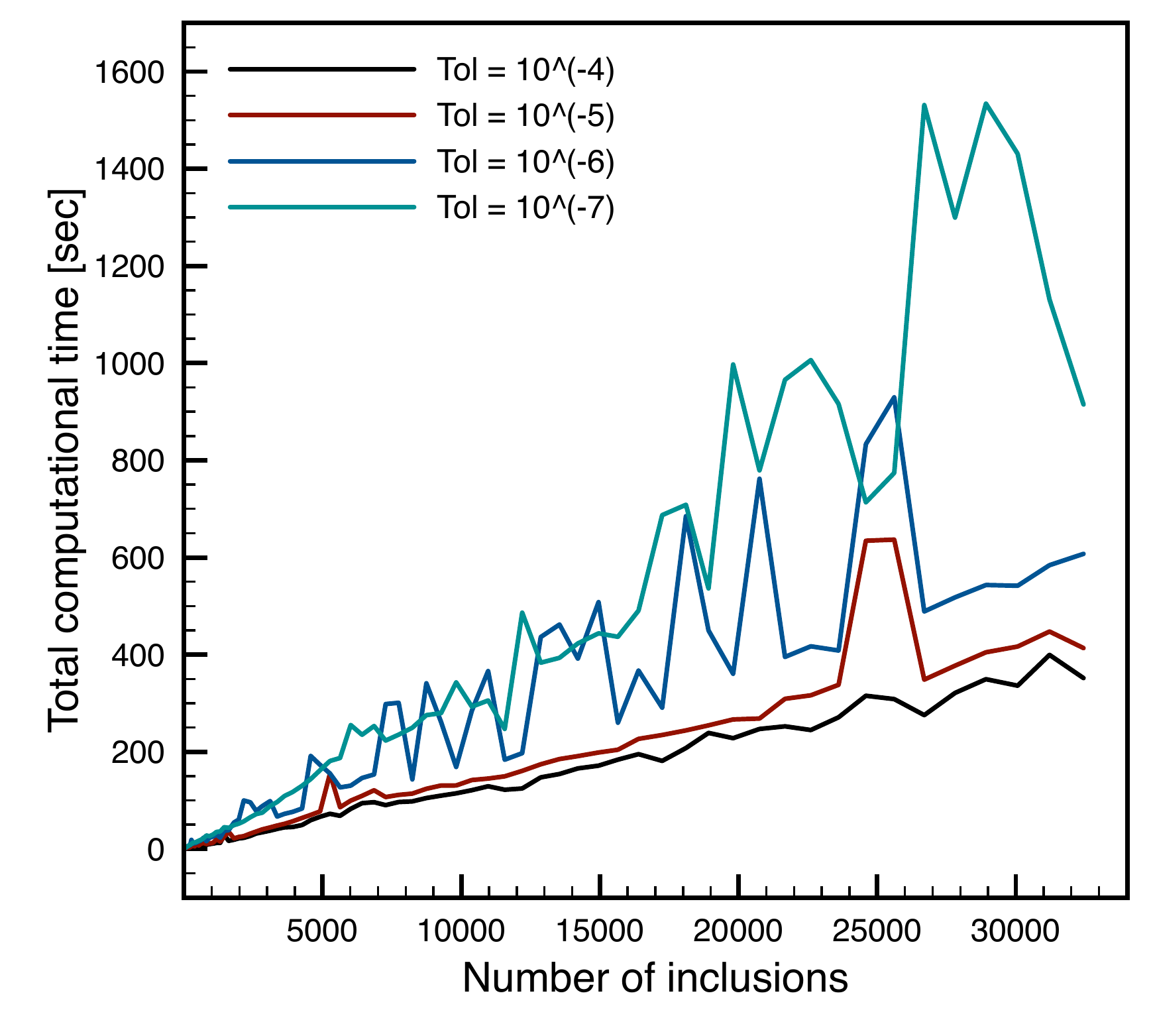}
\includegraphics[scale=0.4]{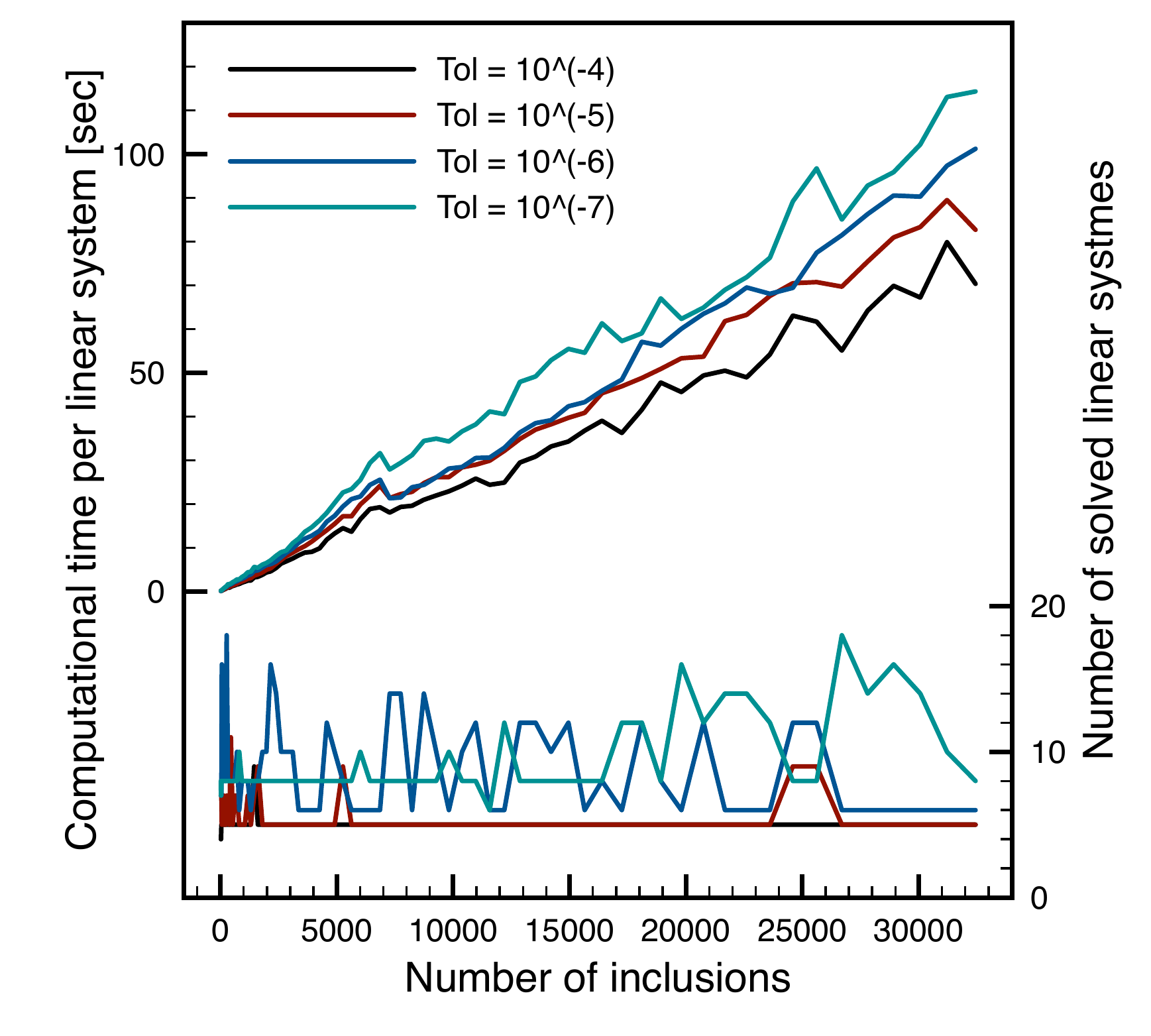}
\caption{[Test case 3] Left: Total computational time to determine the effective diffusion constant. Right: Breakdown into the number of linear systems to be solved and the computational time per linear system.}
\label{fig:T3_bis}
\end{figure}

We conclude this section by presenting some timings. All simulations have been run on a 4~GHz Intel Core~i7 processor, without any parallelization (all computations have been run on a single processor). The method is implemented in Matlab and calls the ScalFMM-library through the MEX-interface. On the left side of Figure~\ref{fig:T3_bis}, we show the wall-clock timings to compute $\ANR{1}$, $\ANR{2}$ and $\ANR{3}$ for different values of the convergence threshold $\eta_{\rm ls}$ and for different numbers of inclusions (the largest system, consisting of 32,442 inclusions, corresponds to $R=20$). We have set $N=1$ for these tests. The threshold $\eta_{\rm opt}$ for the Armijo line seach is chosen 100 times as large, i.e. $\eta_{\rm opt}=100 \, \eta_{\rm ls}$. We observe that the cost increases only linearly with respect to the number $M$ of inclusions, a direct consequence of the use of FMM (without FMM, we expect a cost scaling quadratically with respect to $M$, see Appendix~\ref{ssec:FMM}). 

A breakdown into the number of linear systems to be solved and the typical computational time for solving one such system is shown on the right side of Figure~\ref{fig:T3_bis}. We observe that the time to solve one linear system increases linearly with $M$. However, the number of linear systems to solve is essentially independent of $M$. This shows that the fixed-point procedure (which solves one linear system per iteration) and the optimization algorithm (which solves two linear systems per iterations, due to the use of the Armijo line search) are stable with respect to $M$. Unless very small tolerances $\eta_{\rm opt}$ and $\eta_{\rm ls}$ are chosen, the optimization algorithm converges within a few iterations. 

We also observe that, when the tolerances are set to very small values, the total computational time is not a smooth function of the number of inclusions. This is due to the fact that the number of linear systems to be solved during the optimization procedure varies in a non-smooth manner (in contrast, the average time to solve one linear system increases smoothly with $M$). 

\section{Conclusion}
\label{sec:conc}

In this article, we have proposed a reformulation of the embedded corrector problem in terms of an integral equation of the second type, which, in turn, is discretized using a Galerkin method with real spherical harmonics as basis functions. This numerical method is the building block for the numerical strategy to compute the effective diffusion constants. Thanks to the stability of the formulation with respect to the number of inclusions and the use of an adapted version of the Fast Multipole Method, linear complexity with respect to the number of inclusions is achieved. The various numerical tests that we have performed illustrate the efficiency of the proposed approach.

\appendix

\section{Appendix: Use of the Fast Multipole Method}
\label{ssec:FMM}

We explain here how the Fast Multipole Method (FMM) can be used to efficiently perform the matrix-vector multiplications $\Kv\lambdav$ (see Section~\ref{sec:galerkin}) and $\Kv^\top\sv$ (see Section~\ref{sec:derivative}). Recall that the vectors $\lambdav$ and $\sv$ belong to $\real^{(1+M)(N+1)^2}$ and that $\Kv$ is a square matrix of size $(1+M)(N+1)^2 \times (1+M)(N+1)^2$. A direct computation of $\Kv\lambdav$ or $\Kv^\top\sv$ scales as $M^2$ with respect to the number $M$ of inclusions, due to the nonlocal character of the global potential $\SLO_{\tt G}$. All the other computations of the method can be performed locally on each sphere, and thus in parallel. The computation of $\Kv\lambdav$ or $\Kv^\top\sv$ is thus the bottleneck of the approach in terms of its cost with respect to $M$.

The matrix $\Kv$ is not symmetric, but has a specific structure. Indeed, we have
\[
\Kv = \I - \Lambda \S \Sigma,
\]
where $\I$ is the identity matrix, $\Sigma$ and $\Lambda$ are diagonal matrices given by
\begin{align*}
[\Sigma_{ii}]_{\ell \ell}^{mm} = \frac{\A_0-\A_i}{2\A_0} \, \frac{2\ell+1+\epsilon_i}{4\pi r_i},
\qquad\qquad
[\Lambda_{ii}]_{\ell \ell}^{mm} = \frac{1}{r_i^2},
\end{align*}
and where the matrix $\S$ is given by
\[
[\S_{ij}]_{\ell \ell'}^{mm'} 
= 
r_i^2 \, \frac{4\pi r_j}{2\ell'+1} \sum_{n=1}^{N_g} \omega_n \, \Ylm(s_n) \, (t^{ij}_n)^{f(\ell',\epsilon_j)} \, \Ylpmp\left( s_n^{ij} \right),
\]
where the indices $i$ and $j$ belong to $\IntM= \{1, \ldots, M, \infty\}$. We partition the set $\IntM$ into $\mathcal{M} \cup \{\infty\}$ where $\mathcal{M}:= \{1, \cdots, M\}$ and introduce the following block structure of $\S$:
\[
\S = 
\begin{pmatrix}
\S_{\mathcal M \mathcal M} & \S_{\mathcal M\infty} \\
\S_{\infty \mathcal M} & \S_{\infty\infty} 
\end{pmatrix}.
\]
The matrix-vector products of the matrices $\S_{\mathcal M \infty}$, $\S_{\infty\mathcal M}$ and $\S_{\infty\infty}$ with any vector scale at most linearly with the number $M$ of inclusions. In contrast, the product between the matrix $\S_{\mathcal M \mathcal M}$ and a vector of corresponding size scales {\em quadratically with respect to $M$}, since $\S_{\mathcal M \mathcal M}$ is a non-sparse matrix of size $M(N+1)^2 \times M(N+1)^2$. For any $1 \leq i,j \leq M$, we have
\begin{align*}
[\S_{ij}]_{\ell \ell'}^{mm'} 
&= 
r_i^2 \, \frac{4\pi r_j}{2\ell'+1} \sum_{n=1}^{N_g} \omega_n \, \Ylm(s_n) \, (t^{ij}_n)^{-(\ell'+1)} \, \Ylpmp\left( s_n^{ij} \right).
\end{align*}
We then observe that the matrix $\S_{\mathcal M \mathcal M}$ is almost symmetric, up to numerical integration. Indeed, the symmetric matrix $\widehat \S$ defined by
\begin{align*}
\left[ \widehat \S_{ij} \right]_{\ell \ell'}^{mm'} := \int_{\Gamma_i} \int_{\Gamma_j} \frac{\Ylm(s) \, \Ylpmp(s')}{|s-s'|}ds'ds
\end{align*}
is, up to numerical integration, equal to $\S$: using~\eqref{eq:DefLApp2_Ext}, we have
\begin{align*}
\left[ \widehat \S_{ij} \right]_{\ell \ell'}^{mm'} 
&=
\int_{\Gamma_i} \Ylm(s) \left( \int_{\Gamma_j} \frac{\Ylpmp(s')}{|s-s'|}ds' \right) ds
\\
&=
\int_{\Gamma_i} \Ylm(s) \, (\SL_j \Ylpmp)(s) \, ds
\\
&=
\frac{4\pi r_j}{2\ell'+1} \int_{\Gamma_i} \Ylm(s) \left( \frac{|s-x_j|}{r_j}\right)^{-(\ell'+1)} \Ylpmp\left(\frac{s-x_j}{|s-x_j|}\right) ds
\\
&\approx r_i^2 \, \frac{4\pi r_j}{2\ell'+1} \sum_{n=1}^{N_g} \omega_n \, \Ylm(s_n) \left( \frac{|x_i+r_is_n-x_j|}{r_j}\right)^{-(\ell'+1)} \Ylpmp\left(\frac{x_i+r_is_n-x_j}{|x_i+r_is_n-x_j|}\right)
\\
&=
[\S_{ij}]_{\ell \ell'}^{mm'}.
\end{align*}
As a consequence of this specific structure of $\Kv$, its transpose is approximately given by
\[
\Kv^\top = \I - \Sigma \S^\top \Lambda \approx \I - \Sigma \S \Lambda.
\]
The bottleneck to efficiently compute both $\Kv\lambdav$ and $\Kv^\top\sv$ therefore stems from the matrix $\S_{\mathcal M \mathcal M}$.

\bigskip

For any vector $\kappav \in \real^{M(N+1)^2}$ with coefficients $[\kappa_i]_\ell^m$, the matrix-vector product $\tauv = \S_{\mathcal M \mathcal M} \, \kappav$ is given by
$$
\left[ \tau_i \right]_\ell^m
=
\sum_{j=1}^M \sum_{\ell'=0}^N \sum_{m'=-\ell'}^{\ell'} [\S_{ij}]_{\ell \ell'}^{mm'} [\kappa_j]_{\ell'}^{m'}
=
r_i^2 \sum_{n=1}^{N_g} \omega_n \, V(x_i+r_is_n) \, \Ylm(s_n),
$$
where, for any $x\in \R^3$,
\begin{equation} \label{eq:MultipolarExpansion}
V(x) 
=
\sum_{j=1}^M \sum_{\ell'=0}^N \sum_{m'=-\ell'}^{\ell'} [\Phi_j]_{\ell'}^{m'} \, \frac{1}{|x - x_j|^{\ell'+1}} \, \Ylpmp \left(\frac{x - x_j}{|x - x_j|}\right)
\end{equation}
with 
\[
[\Phi_j]_{\ell'}^{m'} = \frac{4 \pi \, r_j^{\ell'+2}}{2\ell'+1} \, [\kappa_j]_{\ell'}^{m'}.
\]
The function $V$ can be interpreted as the potential generated by $M$ multipoles located at points $x_j$ with momenta $[\Phi_j]_{\ell'}^{m'}$. The quantity $V(x_i+r_is_n)$ is then the value of this potential at the point $x_i+r_is_n$. 

The bottleneck in the evaluation of $\tauv$ thus consists in the efficient evaluation of the potential $V$, created by $M$ multipoles, at a number of points (namely the points $x_i+r_is_n$ for $1 \leq i \leq M$ and $1 \leq n \leq N_g$) that scales proportionally to $M$. Without further approximation, this is clearly an operation that scales quadratically with respect to $M$. However, as we show below, it is possible to resort to the Fast Multipole Method to compute approximations of these quantities in (almost) linear complexity with respect to $M$.

\begin{remark}
As already anticipated, the case of the outer sphere $j=\infty$ needs to be treated separately. The single layer potential inside the sphere $\partial \Bz$ is indeed of the form of a local expansion and contains mode of the form $r^{\ell'}$ instead of $r^{-(\ell'+1)}$. This operation scales linearly (and not quadratically) with respect to the number $M$ of inclusions as only one source is concerned.
\end{remark}

\begin{remark}
Following Remark~\ref{rem:Lij}, we note that the expression~\eqref{eq:Lij} for $i=j$ is a converging approximation of~\eqref{eq:Lii} when $N_g$ increases, and is actually equal to~\eqref{eq:Lii} provided a large enough number of Lebedev points are used. In order not to distinguish the case $i=j$ when using the FMM method, we use only the expression~\eqref{eq:Lij} even if we consider the case $i=j$.
\end{remark}

Standard FMM libraries do not consider arbitrary multipolar expansions as input. The usual case is that only point-charges are considered. In some cases (see e.g.~\cite{librairie-dipole}), dipoles are also treated. Since the degree $N$ of real spherical harmonics that we use is arbitrary, such standard libraries can not be directly used. However, these libraries can be adapted up to some implementation effort. In practice, we have chosen to modify the library ScalFMM~\cite{agullo2014task}.

We assume here that the reader is familiar with the concept of FMM and only comment on the adaptation that is required to use general multipoles as input. An introduction to FMM can be found in~\cite{BeatsonGreengard}. Denoting by $P$ the degree of real spherical harmonics used within the FMM, the following changes need to be made:
\begin{itemize}
\item The P2M-operator needs to be replaced by a M2M-operator, which maps each multipolar expansion of the form~\eqref{eq:MultipolarExpansion} (thus of degree $N$) at the scattered locations $\{x_j\}_{j=1}^M$ in each box to a multipolar expansion of degree $P$ centered at the box. Note that the P2M-operator is a special case of such a M2M-operator with $N=0$.
\item The P2P-operator needs to be replaced by a M2P-operator according to the evaluation of~\eqref{eq:MultipolarExpansion}.
\end{itemize}
The remaining M2L- and L2P-operators remain unchanged.

\begin{remark}
For our numerical tests, we have essentially worked with $N=1$. In that case, FMM librairies able to handle dipoles (such as~\cite{librairie-dipole}) can be used without any further implementation effort.
\end{remark}

\section*{Acknowledgements}

The work of FL is partially supported by ONR under Grant N00014-15-1-2777 and EOARD under grant FA9550-17-1-0294. SX gratefully acknowledges the support from Labex MMCD (Multi-Scale Modelling \& Experimentation of Materials for Sustainable Construction) under contract ANR-11-LABX-0022. 
BS and SX acknowledge the funding from the German Academic Exchange Service (DAAD) from funds of the ``Bundesministeriums f\"ur Bildung und Forschung'' (BMBF) for the project Aa-Par-T (Project-ID 57317909). VE and EC acknowledge the funding from the PICS-CNRS and the PHC PROCOPE 2017 (Project No. 37855ZK). 

\bibliographystyle{abbrv}
\bibliography{biblioPart1,biblio_new}

\end{document}